\RequirePackage{fix-cm}
\documentclass[smallextended]{svjour3}       % onecolumn (second format)
\smartqed  % flush right qed marks, e.g. at end of proof
\usepackage{graphicx}
\usepackage[misc]{ifsym}
\usepackage{multirow}
\usepackage{amsmath,amssymb}
\usepackage{float}
\usepackage{url}
\usepackage[ruled]{algorithm2e}
\usepackage[top=1.0in, bottom=1.0in, left=0.9in, right=0.9in]{geometry}
\usepackage{caption}
\usepackage{algpseudocode}
\usepackage[colorlinks,linkcolor=blue,citecolor=blue]{hyperref}
\if@onethmnum
\newtheorem{assumption}[theorem]{Assumption}
\else
\newtheorem{assumption}[theorem]{Assumption}
\fi
\usepackage{bm}

\begin{document}

\markboth{X.~L.~Song, B.~Chen and B.~Yu}{Error Estimates for Sparse Optimal Control Problems}

\title{Error Estimates for Sparse Optimal Control Problems by Piecewise Linear Finite Element Approximation}

\author{Xiaoliang Song
\thanks{School of Mathematical Sciences,
Dalian University of Technology, Dalian, Liaoning, China\\ Email: songxiaoliang@mail.dlut.edu.cn}\and
Bo Chen
\thanks{Department of Mathematics, National University of Singapore, 10 Lower Kent Ridge Road, Singapore\\ Email: chenbo@u.nus.edu} \and
Bo Yu
\thanks{School of Mathematical Sciences,
Dalian University of Technology, Dalian, Liaoning, China\\ Email: yubo@dlut.edu.cn}
}

\maketitle

\begin{abstract}
Optimization problems with $L^1$-control cost functional subject to an elliptic partial differential equation (PDE) are considered. However, different from the finite dimensional $l^1$-regularization optimization, the resulting discretized $L^1$-norm does not have a decoupled form when the standard piecewise linear finite element is employed to discretize the continuous problem. A common approach to overcome this difficulty is employing a nodal quadrature formula to approximately discretize the $L^1$-norm. It is inevitable that this technique will incur an additional error. Different from the traditional approach, a duality-based approach
and an accelerated block coordinate descent (ABCD) method is introduced to solve this type of problem via its dual. Based on the discretized dual problem, a new discretized scheme for the $L^1$-norm is presented.  Compared new discretized scheme for $L^1$-norm with the nodal quadrature formula, the advantages of our new discretized scheme can be demonstrated in terms of the approximation order. More importantly, finite element error estimates results for the primal problem with the new discretized scheme for the $L^1$-norm are provided, which confirm that this approximation scheme will not change the order of error estimates.
\end{abstract}

%\begin{classification}
%49N05 65N30 49M25 68W15
%\end{classification}

\begin{keywords}{finite element method, ABCD method, approximate discretization, error estimates.}
\end{keywords}

\section{Introduction}\label{sec:1}
In this paper, we study the following linear-quadratic elliptic PDE-constrained optimal control problem with $L^1$-control cost and piecewise box constraints on the control
\begin{equation}\label{eqn:orginal problems}
           \qquad \left\{ \begin{aligned}
        &\min \limits_{(y,u)\in Y\times U}^{}\ \ J(y,u)=\frac{1}{2}\|y-y_d\|_{L^2(\Omega)}^{2}+\frac{\alpha}{2}\|u\|_{L^2(\Omega)}^{2}+\beta\|u\|_{L^1(\Omega)} \\
        &\qquad{\rm s.t.}\qquad Ly=u+y_r\ \ \mathrm{in}\  \Omega, \\
         &\qquad \qquad \qquad y=0\qquad\quad  \mathrm{on}\ \partial\Omega,\\
         &\qquad \qquad\qquad u\in U_{ad}=\{v(x)|a\leq v(x)\leq b,\ {\rm a.e. }\  \mathrm{on}\ \Omega\}\subseteq U,
                          \end{aligned} \right.\tag{$\mathrm{P}$}
 \end{equation}
where $Y:=H_0^1(\Omega)$, $U:=L^2(\Omega)$, $\Omega\subseteq \mathbb{R}^n$ ($n=2$ or $3$) is a convex, open and bounded domain with $C^{1,1}$- or polygonal boundary $\Gamma$, the desired state $y_d\in L^2(\Omega)$ and the source term $y_r \in L^2(\Omega)$ are given; and $a\leq0\leq b$, $\alpha$, $\beta>0$. Moreover, the operator $L$ is a second-order linear elliptic differential operator. It is well-known that $L^1$-norm could lead to sparse optimal control, i.e. the optimal control with small support. Such an optimal control problem (\ref{eqn:orginal problems}) plays an important role for the placement of control devices \cite{Stadler}. In some cases, it is difficult or undesirable to place control devices all over the control domain and one hopes to localize controllers in small and effective regions, the $L^1$-solution gives information about the optimal location of the control devices.

Let us comment on known results on a-priori analysis of control constrained sparse optimal control problems. For the study of optimal control problems with sparsity promoting terms, as far as we know, the first paper devoted to this study is published by Stadler \cite{Stadler}, in which structural properties of the control variables were analyzed in the case of the linear-quadratic elliptic optimal control problem. In 2011, a priori and a posteriori error estimates were first given by G. Wachsmuth and D. Wachsmuth in \cite{WaWa} for piecewise linear control discretizations, in which they proved the following result
\begin{equation*}
  \|u^*-u^*_h\|_{L^2(\Omega)}\leq C(\alpha^{-1}h+\alpha^{-3/2}h^2).
\end{equation*}
However, from an algorithmic point of view, using the piecewise linear finite elements with nodal basis functions $\{\phi_i(x)\}$ to approximate the control variable $u$, the resulting discretized $L^1$-norm
\begin{eqnarray*}\label{equ:discrete norm}
\|u_h\|_{L^1(\Omega_h)}&=&\int_{\Omega_h}\big{|}\sum\limits_{i=1}^{N_h}u_i\phi_i(x)\big{|}\mathrm{d}x,
\end{eqnarray*}
lead to its subgradient $\nu_h\in \partial\|u_h\|_{L^1(\Omega_h)}$ not to be expressed by $\{\phi_i(x)\}$ since $\nu_h$ may have jumps along lines $u_h=0$ which are not grid lines. In addition, the discretized $L^1$-norm does not have a decoupled form with respect to the coefficients $\{u_i\}$.
Thus, directly solving the corresponding discretized problem will cause many difficulties in numerical calculation.
Hence, the authors in \cite{WaWa} introduced an alternative discretization of the $L^1$-norm which relies on a nodal quadrature formula
\begin{eqnarray}\label{approx discretization11}
\|u_h\|_{L^{1}_h(\Omega_h)}&:=&\sum_{i=1}^{n}|u_i|\int_{\Omega_h}\phi_i(x)\mathrm{d}x.
\end{eqnarray}
About the approximate $L^1$-norm, based on the error estimates of the nodal interpolation operator, it is easy to show that
\begin{equation}\label{equ:approximate discrete norm error}
  0\leq\|u_h\|_{L^{1}_h(\Omega_h)}-\|u_h\|_{L^{1}(\Omega_h)}=O(h).
\end{equation}
Obviously, this quadrature incurs an additional error. However, the authors \cite{WaWa} proved that this approximation does not change the order of error estimates and they showed that
\begin{equation*}
  \|u^*-u^*_h\|_{L^2(\Omega)}\leq C(\alpha^{-1}h+\alpha^{-3/2}h^2).
\end{equation*}

In a sequence of papers \cite{CaHerWa1,CaHerWa2}, for the non-convex case governed by a semilinear elliptic equation, Casas et al. proved second-order necessary and sufficient optimality conditions. Using the second-order sufficient optimality conditions, the authors provided error estimates of order $h$ w.r.t. the $L^\infty$ norm for three different choices of the control discretization (including the piecewise constant, piecewise linear control discretization and the variational control discretization ). It should be pointed out that, for the piecewise linear control discretization case,
a similar approximation technique to the one introduced by G. Wachsmuth and D. Wachsmuth is also used for the discretizations of the $L^2$-norm and $L^1$-norm of the control.

To numerically solve the problem (\ref{eqn:orginal problems}), there are two possible ways.
One is called \emph{First discretize, then optimize}. More specifically, one may first discretize the continuous problem by using the finite element method, which results in a finite dimensional optimization problem. Then, the corresponding finite dimensional optimization problem can be solved numerically with the help of a suitable algorithm. Instead of applying discretized concepts to the continuous problem directly, another approach is first applying an algorithm on the continuous level or computing the infinite dimensional optimality system, and then discretizing the related subproblems appeared in the algorithm or the optimality system by using the finite element method. This approach is called \emph{First optimize, then discretize}.  There are different opinions regarding which route to take (see Collis and Heinkenschloss \cite{CollisHeink} for a discussion). Independently of where discretization is located, the resulting finite dimensional equations are quite large. Thus, both of these cases require us to consider proposing an efficient algorithm.

Let us mention some existing numerical methods for solving the optimal control problem (\ref{eqn:orginal problems}). Since the problem (\ref{eqn:orginal problems}) is nonsmooth, thus applying semismooth Newton (SSN) method is used to be a priority in consideration of their locally superlinear convergence. SSN method in function space is proved to have the locally superlinear convergence (see \cite{Ulbrich1,Ulbrich2,HiPiUl} for more details). Furthermore, mesh-independence results for SSN method were established in \cite{meshindependent}.

Although employing the SSN method can derive the solution with high precision, it is generally known that the total error of utilizing numerical methods to solve PDE constrained problem consists of two parts: discretization error and the iteration error resulted from algorithm to solve the discretized problem. Since the error order of piecewise linear finite element method is $\mathcal{O}(h)$ which accounts for the main part, algorithms of high precision do not reduce the order of the total error but waste computations. Thus, taking the precision of discretization error into account, employing fast and efficient first-order algorithms with the aim of solving discretized problems to moderate accuracy is sufficient. As one may know, for finite dimensional large scale optimization problems, some efficient first-order algorithms, such as iterative soft thresholding algorithms (ISTA) \cite{Blumen}, accelerated proximal gradient (APG)-based method \cite{inexactAPG,Beck,Toh}, alternating direction method of multipliers (ADMM) \cite{Fazel,SunToh1,SunToh2}, etc, have become the state of the art algorithms. Motivated by the success of these finite dimensional optimization algorithms, an ADMM \cite{iwADMM} and an APG \cite{FIP} method are proposed in function space to solve the sparse optimal control problems.

As far as we know, most of the aforementioned papers are devoted to solve the primal problem (\ref{eqn:orginal problems}). However, as mentioned above, from the perspective of actual numerical implementation, directly solving the primal problem is difficult, since the discretized $L^1$-norm does not have a decoupled form
when the primal problem (\ref{eqn:orginal problems}) is discretized by the piecewise linear finite element. Thus the same technique as (\ref{approx discretization11}) should be used, which will inevitably cause additional error. Alternatively, instead of solving the primal problem, in \cite{mABCDSOPT}, Song et al. considered using the duality-based approach for (\ref{eqn:orginal problems}) and solving the dual problem. Taking advantage of the structure of the dual problem, the authors employed a majorized accelerated block coordinate descent (mABCD) method to solve the dual problem. Specifically, combining an inexact 2-block majorized ABCD \cite[Chapter 3]{CuiYing} and the recent advances in the inexact symmetric Gauss-Seidel (sGS) decomposition technique developed in \cite{SunToh2,SunToh3}, Song et al. proposed a sGS based majorized ABCD method (called sGS-mABCD) to solve the dual problem.

Owing to the important convergence results of mABCD method, in \cite{mABCDSOPT}, the sGS-mABCD algorithm builds a sequence of iterations
$\{{\bf{z}}^k\}:=\{({\bm{\lambda}}^k,{\bf p}^k,{\bm \mu}^k)\}$ for which $\Phi_h({\bf{z}}^k)- \Phi_h({\bf z}^*)=O(1/k^2)$, where $\Phi_h$ is the dual objective function. Based on the second order growth condition of $\Phi_h$, Song et al. in \cite{meshindependencemABCD} also showed that
\begin{equation*}
  {\rm dist}({\bf z}^k,(\partial\Phi_h)^{-1}(0))=O(1/k).
\end{equation*}
$\bf\lambda$
More importantly, in \cite{meshindependencemABCD}, the authors also gave two types of  mesh-independence for mABCD method, which assert that asymptotically the infinite dimensional mABCD method and finite
dimensional discretizations have the same convergence property, and the iterations of mABCD method remain
nearly constant as the discretization is refined.

Although the convergence behavior and the iteration complexity of the dual problem have been shown, our ultimate goal is looking for the optimal control solution. Thus, it is necessary to analyze the primal problem. As shown in Section \ref{sec:4.2}, the primal problem of the discretized dual problem is also an approximate discretization of problem (\ref{eqn:orginal problems}), in which the $L^1$-norm is approximated by
\begin{equation}\label{approx discretization21}
\|u_h\|_{\widetilde{L}^1_h(\Omega)}=\sum_{j=1}^{N_h}|\int_{\Omega_h}\sum\limits_{i=1}^{N_h}u_i\phi_i\phi_j(x)\mathrm{d}x|=\|M_h{\bf u}\|_1,
\end{equation}
where $M_h$ is the mass matrix and ${\bf u}=(u_1,u_2,...,u_{N_h})$. Since $\|u_h\|_{\widetilde{L}^1_h(\Omega)}$ can be regarded as an approximation of $\int_{\Omega_h}|\sum_{i=1}^{N_h}{u_i\phi_i(x)}|~\mathrm{d}x$, it is necessarily required to analyse the relationship between them. In this paper, we can show the following result %(see Proposition \ref{eqn:another martix properties}):
\begin{equation}\label{new approx L1 norm error}
0\leq\|u_h\|_{L^{1}(\Omega_h)}-\|u_h\|_{\widetilde{L}^{1}_h(\Omega_h)}=O(h^2).
\end{equation}
More importantly, another key issue should be considered is how measures of the solution
accuracy by using discretized form (\ref{approx discretization21}) vary with the level of discretized approximation. Such questions come under the category of the finite element error estimates. In this paper, we will explain the reasonability of employing discretized form (\ref{approx discretization21}) and give our main important error estimates results (see Theorem \ref{theorem:error2} and Corollary \ref{corollary:error1}):
\begin{equation*}
  \|u^*-u^*_h\|_{L^2(\Omega)}\leq C(\alpha^{-1}h+\alpha^{-3/2}h^2).
\end{equation*}

Obviously, employing discretized form (\ref{approx discretization21}) will also not change the order of error estimates. Thus, compared (\ref{approx discretization11}) with (\ref{approx discretization21}), it is obvious that utilizing $\|u_h\|_{\widetilde{L}_h^1(\Omega_h)}$ to approximate $\|u_h\|_{L^1(\Omega_h)}$ is better than using $\|u_h\|_{L_h^1(\Omega_h)}$ in term of the approximation order. Hence, equivalently solving the dual problem with the discretized form (\ref{approx discretization21}) is superior to directly solving the primal problem with the discretized form (\ref{approx discretization11}). Actually, in \cite{mABCDSOPT}, from their numerical results, the authors have already shown that solving the dual problem could get better error results than that from solving the primal problem.

The remainder of the paper is organized as follows. In Section \ref{sec:2}, the first-order optimality conditions for problem (\ref{eqn:orginal problems}) are derived. In Section \ref{sec:3}, piecewise linear finite element discretization and an approximate discretization approach are introduced. In Section \ref{sec:4}, we give a brief sketch of the symmetric Gauss-Seidel based majorized ABCD (sGS-mABCD) method for the dual problem, and show some convergence results. More importantly, based on the discretized dual problem, a new approximate discretization primal problem is presented. In Section \ref{sec:5}, some error estimates results are proved for the new approximate discretization primal problem. Finally, we conclude our paper in Section \ref{sec:6}.

\section{First-order optimality condition}
\label{sec:2}
In this section, we will derive the first-order optimality conditions. Firstly, let us suppose the elliptic PDEs involved in (\ref{eqn:orginal problems}) which are of the form
\begin{equation}\label{eqn:state equations}
\begin{aligned}
 Ly&=u+y_r \quad \mathrm{in}\  \Omega, \\
 y&=0  \qquad \quad\ \mathrm{on}\ \partial\Omega,
\end{aligned}
\end{equation}
satisfy the following assumption.
\begin{assumption}\label{equ:assumption:1}
The linear second-order differential operator $L$ is defined by
 \begin{equation}\label{operator A}
   (Ly)(x):=-\sum \limits^{n}_{i,j=1}\partial_{x_{j}}(a_{ij}(x)y_{x_{i}})+c_0(x)y(x),
 \end{equation}
where functions $a_{ij}(x), c_0(x)\in L^{\infty}(\Omega)$, $c_0\geq0$,
and it is uniformly elliptic, i.e. $a_{ij}(x)=a_{ji}(x)$ and there is a constant $\theta>0$ such that
\begin{equation}\label{equ:operator A coercivity}
  \sum \limits^{n}_{i,j=1}a_{ij}(x)\xi_i\xi_j\geq\theta\|\xi\|^2 \quad \mathrm{for\ a.a.}\ x\in \Omega\  \mathrm{and}\  \forall \xi \in \mathbb{R}^n.
\end{equation}
\end{assumption}

The weak formulation of (\ref{eqn:state equations}) is given by
\begin{equation}\label{eqn:weak form}
  \mathrm{Find}\ y\in H_0^1(\Omega):\ a(y,v)=(u+y_r,v)_{L^2(\Omega)},\quad \forall v \in H_0^1(\Omega),
\end{equation}
with the bilinear form
\begin{equation}\label{eqn:bilinear form}
  a(y,v)=\int_{\Omega}(\sum \limits^{n}_{i,j=1}a_{ji}y_{x_{i}}v_{x_{i}}+c_0yv)\mathrm{d}x.
\end{equation}
%or in short $ Ay=B(u+y_r)$, where $A\in \mathcal{L}(Y,Y^*)$ is the operator induced by the bilinear form $a$, i.e., $Ay=a(y,\cdot)$ and $B\in \mathcal{L}(U,Y^*)$ is defined by $Bu=(u,\cdot)_{L^2(\Omega)}$. Since the bilinear form $a(\cdot,\cdot)$ is symmetric and $U,Y$ are Hilbert spaces, we have $A^*\in\mathcal{L}(Y,Y^*)=A$, and $B^*\in\mathcal{L}(Y,U)$ with $B^*v=v, \forall v\in Y$.

\begin{remark}\label{more general case}
Although we assume that the Dirichlet boundary condition $y=0$ holds, it should be noted that the
assumption is not a restriction and our considerations can also carry over to the more general boundary conditions of Robin type
\begin{equation*}
  \frac{\partial y}{\partial \nu}+\gamma y=g \quad {\rm on}\  \partial\Omega,
\end{equation*}
where $g\in L^2(\partial\Omega)$ is given and $\gamma\in L^{\infty}(\partial\Omega)$ is nonnegative coefficient. Furthermore, it is assumed that the control satisfies $a\leq u \leq b$, where $a$ and $b$ have opposite signs. First, we should emphasize that this condition is required in practice, e.g., the placement of control devices. In addition, please also note, that this condition is not a restriction from the point of view of the algorithm.
If one has, e.g., $a>0$ on $\Omega$, the $L^1$-norm in $U_{ad}$ is in fact a linear
function, and thus the problem can also be handled in an analogous way.
\end{remark}

Then, we analyze the existence and uniqueness of global solution to problem (\ref{eqn:orginal problems}). Utilizing the Lax-Milgram lemma, we have the following proposition.
\begin{proposition}{\rm\textbf{\cite [Theorem B.4] {KiSt}}}\label{equ:weak formulation}
%{\rm(Existence and Uniqueness for the Dirichlet Problem)}
Under Assumption {\ref{equ:assumption:1}}, the bilinear form $a(\cdot,\cdot)$ in {\rm (\ref{eqn:bilinear form})} is bounded and $V$-coercive for $V=H^1_0(\Omega)$. In particular, for every $u \in L^2(\Omega)$ and $y_r\in L^2(\Omega)$, {\rm (\ref{eqn:state equations})} has a unique weak solution $y\in H^1_0(\Omega)$ given by {\rm (\ref{eqn:weak form})}. Furthermore,
\begin{equation}\label{equ:control estimats state}
  \|y\|_{H^1}\leq C (\|u\|_{L^2(\Omega)}+\|y_r\|_{L^2(\Omega)}),
\end{equation}
for a constant $C$ depending only on $a_{ij}$, $c_0$ and $\Omega$.
\end{proposition}

By Proposition \ref{equ:weak formulation}, the weak formulation (\ref{eqn:weak form}) can be rewritten as follows:
\begin{equation*}
   Ay=B(u+y_r),
\end{equation*}
where $A\in \mathcal{L}(Y,Y^*)$ is the operator induced by the bilinear form $a$, i.e. $Ay=a(y,\cdot)$ and $B\in \mathcal{L}(U,Y^*)$ is defined by $Bu=(u,\cdot)_{L^2(\Omega)}$. Since the bilinear form $a(\cdot,\cdot)$ is symmetric and $U,Y$ are Hilbert spaces, we have $A^*\in\mathcal{L}(Y,Y^*)=A$ and $B^*\in\mathcal{L}(Y,U)$ with $B^*v=v$ for any $v\in Y$.

Furthermore, for the convenience of later error estimates, we introduce the solution operator $\mathcal{S}$: $H^{-1}(\Omega)\rightarrow H^1_0(\Omega)$ with $y(u):=\mathcal{S}(u+y_r)$, which is called the control-to-state mapping and is a continuous linear injective operator. Since $H_0^1(\Omega)$ is a Hilbert space, the adjoint operator $\mathcal{S^*}$: $H^{-1}(\Omega)\rightarrow H_0^1(\Omega)$ is also a continuous linear operator.

From the strong convexity and lower semicontinuity of the objective functional $J(y,u)$ of (\ref{eqn:orginal problems}), it is easy to establish the existence and uniqueness of the solution to (\ref{eqn:orginal problems}). The optimal solution can be characterized by the following Karush-Kuhn-Tucker (KKT) conditions:

\begin{theorem}[{\rm First-Order Optimality Condition}]\label{First-order optimality condition}
Under Assumption \ref{equ:assumption:1}, ($y^*$, $u^*$) is the optimal solution of {\rm(\ref{eqn:orginal problems})}, if and only if there exists an adjoint state $p^*\in H_0^1(\Omega)$, such that the following conditions hold in the weak sense
\begin{subequations}\label{eqn:KKT}
\begin{eqnarray}
&&\begin{aligned}\label{eqn1:KKT}
%        Ay^*&=u^*+y_r \quad\mathrm{in}\ \Omega,\\
%        y^*&=0 \ ~\quad\qquad\mathrm{on}\ \partial\Omega,
        y^*=\mathcal{S}(u^*+y_r),
        \end{aligned}  \\
&& \begin{aligned}\label{eqn2:KKT}
%        Ap^*&=y_d-y^*&  &\mathrm{in}\ \Omega,\\
%        p^*&=0& &\mathrm{on}\ \partial\Omega,
        p^*&=\mathcal{S}^*(y_d-y^*),
        \end{aligned}\\
&&{\left\langle\alpha u^*-p^*,u-u^*\right\rangle_{L^2(\Omega)}+\beta(\|u\|_{L^1(\Omega)}-\|u^*\|_{L^1(\Omega)})}\geq0,\quad \forall u \in U_{ad}.\label{eqn3:KKT}
\end{eqnarray}
\end{subequations}
\end{theorem}

\begin{remark}\label{relation}
It is easy to obtain that the variational inequality (\ref{eqn3:KKT}) can be equivalently rewritten as the following nonsmooth equation:
\begin{equation}\label{eqn3.1:KKT}
  u^*={\rm\Pi}_{U_{ad}}\left(\frac{1}{\alpha}{\rm{soft}}\left(p^*,\beta\right)\right),
\end{equation}
where
\begin{equation*}
  \begin{aligned}
 {\rm\Pi}_{U_{ad}}(v(x))&:=\max\{a,\min\{v(x),b\}\}, \\
 {\rm soft}(v(x),\beta)&:={\rm{sgn}}(v(x))\circ\max(|v(x)|-\beta,0).
  \end{aligned}
\end{equation*}
From (\ref{eqn3.1:KKT}), an obvious fact should be pointed out that $|p|<\beta$ implies $u=0$, which also explains that the $L^1$-norm can induce the sparsity property of $u$. Moreover, since $p\in H^1_0(\Omega)$, (\ref{eqn3:KKT}) also implies $u\in H^1(\Omega)$. Figure \ref{fig:the relation u and p} shows the relationship between $u$ and $p$.
\begin{figure}[h]
\centering
\includegraphics[width=0.40\textwidth]{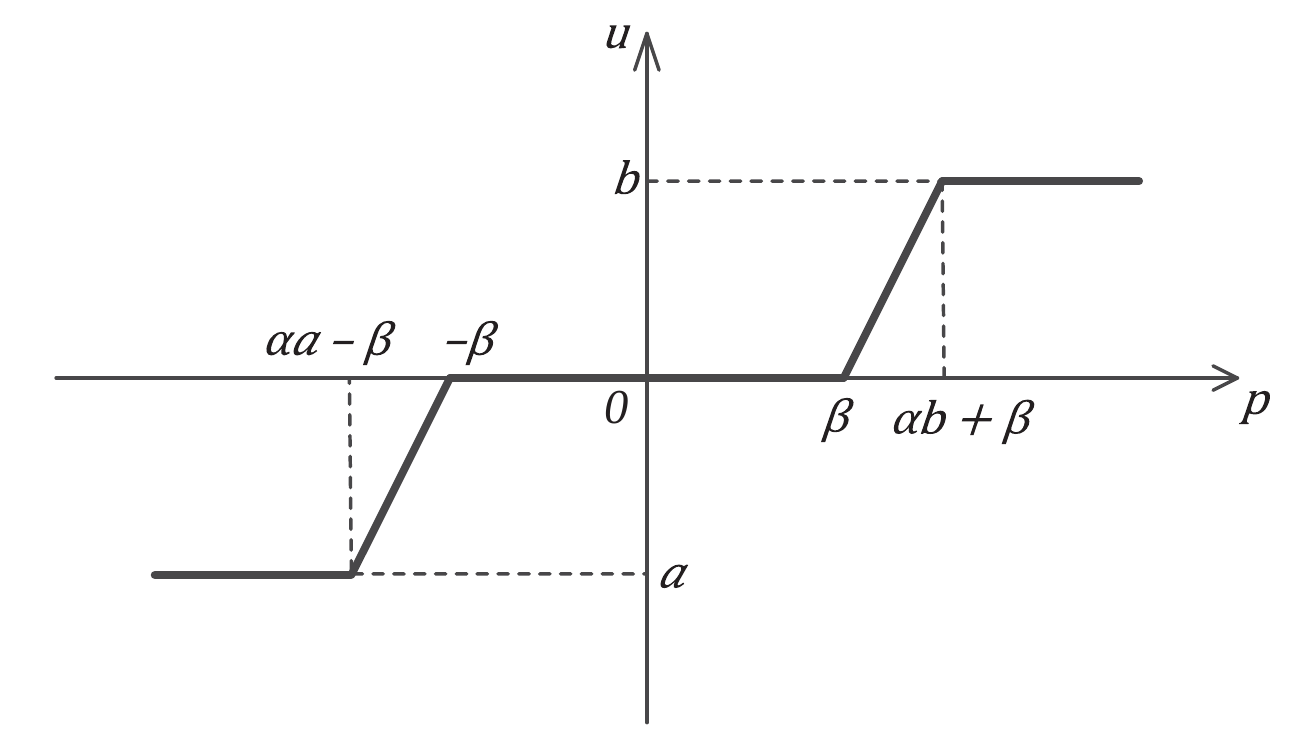}
\caption{ The relationship between $u$ and $p$}\label{fig:the relation u and p}
\end{figure}
\end{remark}
It is obvious that if $\beta$ is sufficiently large, the optimal control would be $u^*_{\beta}=0$. Then we have the following lemma.
\begin{lemma}
  If $\beta\geq\beta_0:=\|A^{-*}(y_d-A^{-1}By_r)\|_{L^\infty{(\Omega)}}$, then the unique solution of problem {\rm(\ref{eqn:orginal problems})} is $(y^*_{\beta},u^*_{\beta})=(A^{-1}By_r,0)$.
\proof By Proposition {\rm\ref{equ:weak formulation}}, we know that the operator $A$ has a bounded inverse, then we can give the following reduced objective function $\widehat{J}$:
\begin{equation}\label{reduced objective function }
  \hat{J}(u):=J(A^{-1}B(u+y_r),u)=\frac{1}{2}\|A^{-1}Bu+A^{-1}By_r-y_d\|_{L^2(\Omega)}^{2}+\frac{\alpha}{2}\|u\|_{L^2(\Omega)}^{2}+\beta\|u\|_{L^1(\Omega)}.
\end{equation}
For any $u\in L^2(\Omega)$, we have
\begin{equation*}
  \begin{aligned}
    \hat{J}(u)-\hat{J}(0) &=\frac{1}{2}\|A^{-1}Bu\|_{L^2(\Omega)}^{2}-\langle y_d-A^{-1}By_r, A^{-1}Bu\rangle_{L^2(\Omega)}+\frac{\alpha}{2}\|u\|_{L^2(\Omega)}^{2}+\beta\|u\|_{L^1(\Omega)} \\
     %&= \frac{1}{2}\|A^{-1}Bu\|_{L^2(\Omega)}^{2}-\langle A^{-*}(y_d-A^{-1}By_r), u\rangle_{L^2(\Omega)}+\frac{\alpha}{2}\|u\|_{L^2(\Omega)}^{2}+\beta\|u\|_{L^1(\Omega)}\\
     &\geq\frac{1}{2}\|A^{-1}Bu\|_{L^2(\Omega)}^{2}-\|u\|_{L^1(\Omega)}\|A^{-*}(y_d-A^{-1}By_r)\|_{L^\infty{(\Omega)}} +\frac{\alpha}{2}\|u\|_{L^2(\Omega)}^{2}+\beta\|u\|_{L^1(\Omega)}\\
     &=\frac{1}{2}\|A^{-1}Bu\|_{L^2(\Omega)}^{2} +\frac{\alpha}{2}\|u\|_{L^2(\Omega)}^{2}+(\beta-\beta_0)\|u\|_{L^1(\Omega)}
  \end{aligned}
\end{equation*}
Obviously, if $\beta\geq\beta_0$, the latter expression is nonnegative. Thus, for $\beta\geq\beta_0$, $\hat{J}(u)\geq\hat{J}(0)$ for all $u\in U_{ad}$, which proves that the optimal control is $u^*_{\beta}=0$ and the corresponding state $y^*_{\beta}=A^{-1}By_r$.
\end{lemma}

\section{Finite element discretization}
\label{sec:3}
To numerically solve problem (\ref{eqn:orginal problems}), we consider the finite element method, in which the state $y$ and the control $u$ are both discretized by the piecewise linear, globally continuous finite elements.

To this aim, let us fix the assumptions on the discretization by finite elements. We first consider a family of regular and quasi-uniform triangulations $\{\mathcal{T}_h\}_{h>0}$ of $\bar{\Omega}$. For each cell $T\in \mathcal{T}_h$, let us define the diameter of the set $T$ by $\rho_{T}:={\rm diam}\ T$ and define $\sigma_{T}$ to be the diameter of the largest ball contained in $T$. The mesh size of the grid is defined by $h=\max_{T\in \mathcal{T}_h}\rho_{T}$. We suppose that the following regularity assumptions on the triangulation are satisfied which are standard in the context of error estimates.

\begin{assumption}[regular and quasi-uniform triangulations]\label{regular and quasi-uniform triangulations}
There exist two positive constants $\kappa$ and $\tau$ such that
   \begin{equation*}
   \frac{\rho_{T}}{\sigma_{T}}\leq \kappa,\quad \frac{h}{\rho_{T}}\leq \tau,
 \end{equation*}
hold for all $T\in \mathcal{T}_h$ and all $h>0$. Moreover, let us define $\bar{\Omega}_h=\bigcup_{T\in \mathcal{T}_h}T$, and let ${\Omega}_h \subset\Omega$ and $\Gamma_h$ denote its interior and its boundary, respectively. In the case that $\Omega$ is a convex polyhedral domain, we have $\Omega=\Omega_h$. In the case $\Omega$ has a $C^{1,1}$- boundary $\Gamma$, we assume that $\bar{\Omega}_h$ is convex and all boundary vertices of $\bar{\Omega}_h$ are contained in $\Gamma$, such that
\begin{equation*}
  |\Omega\backslash {\Omega}_h|\leq c h^2,
\end{equation*}
where $|\cdot|$ denotes the measure of the set and $c>0$ is a constant.
\end{assumption}

\subsection{\textbf {Piecewise linear finite elements discretization}}
\label{sec:3.1}
On account of the homogeneous boundary condition of the state equation, we use
\begin{equation}\label{state discretized space}
  Y_h =\left\{y_h\in C(\bar{\Omega})~\big{|}~y_{h|T}\in \mathcal{P}_1~ {\rm{for\ all}}~ T\in \mathcal{T}_h~ \mathrm{and}~ y_h=0~ \mathrm{in } ~\bar{\Omega}\backslash {\Omega}_h\right\}
\end{equation}
as the discretized state space, where $\mathcal{P}_1$ denotes the space of polynomials of degree less than or equal to $1$. For a given source term $y_r$ and right-hand side $u\in L^2(\Omega)$, we denote by $y_h(u)$ the approximated state associated with $u$, which is the unique solution for the following discretized weak formulation: %of the state equation (\ref{eqn:state equations}):
\begin{equation}\label{eqn:discrete weak solution}
 \int_{\Omega_h}\left(\sum \limits^{n}_{i,j=1}a_{ij}{y_h}_{x_{i}}{v_h}_{x_{j}}+c_0y_hv_h\right)\mathrm{d}x=\int_{\Omega_h}(u+y_r)v_h{\rm{d}}x \quad  \forall v_h\in Y_h.
\end{equation}
Moreover, $y_h(u)$ can also be expressed by $y_h(u)={\mathcal{S}}_{h}(u+y_r)$, in which  ${\mathcal{S}}_{h}$ is a discretized vision of $\mathcal{S}$ and an injective, selfadjoint operator. The following error estimates are well-known.
\begin{lemma}{\rm\textbf{\cite[Theorem 4.4.6]{Ciarlet}}}\label{eqn:lemma1}
For a given $u\in L^2(\Omega)$, let $y$ and $y_h(u)$ be the unique solution of {\rm(\ref{eqn:weak form})} and {\rm(\ref{eqn:discrete weak solution})}, respectively. Then there exists a constant $c_1>0$ independent of $h$, $u$ and $y_r$ such that

\begin{equation}\label{estimates1}
  \|y-y_h(u)\|_{L^2(\Omega)}+h\|\nabla y-\nabla y_h(u)\|_{L^2(\Omega)}\leq c_1h^2(\|u\|_{L^2(\Omega)}+\|y_r\|_{L^2(\Omega)}).
\end{equation}
In particular, this implies $\|\mathcal{S}-\mathcal{S}_h\|_{L^2\rightarrow L^2}\leq c_1h^2$ and $\|\mathcal{S}-\mathcal{S}_h\|_{L^2\rightarrow H^1}\leq c_1h$.
\end{lemma}

As mentioned above, we also use the same discretized space to discretize control $u$, thus we define
\begin{equation}\label{control discretized space}
   U_h =\left\{u_h\in C(\bar{\Omega})~\big{|}~u_{h|T}\in \mathcal{P}_1~ {\rm{for\ all}}~ T\in \mathcal{T}_h~ \mathrm{and}~ u_h=0~ \mathrm{in } ~\bar{\Omega}\backslash{\Omega}_h\right\}.
\end{equation}
For a given regular and quasi-uniform triangulation $\mathcal{T}_h$ with nodes $\{x_i\}_{i=1}^{N_h}$, let $\{\phi_i\} _{i=1}^{N_h}$ be a set of nodal basis functions, which span $Y_h$ as well as $U_h$ and satisfy the following properties
\begin{eqnarray}\label{basic functions properties}
  &&\phi_i \geq 0, \quad
  \|\phi_i\|_{\infty} = 1 \quad \forall i=1,2,...,N_h,
 \quad \sum\limits_{i=1}^{N_h}\phi_i(x)=1.
\end{eqnarray}
The elements $y_h\in Y_h$ and $u_h\in U_h$can be represented in the following forms, respectively,
\begin{equation*}
  y_h=\sum \limits_{i=1}^{N_h}y_i\phi_i(x),\quad u_h=\sum \limits_{i=1}^{N_h}u_i\phi_i(x),
\end{equation*}
where $y_h(x_i)=y_i$ and $u_h(x_i)=u_i$. Moreover, we use ${\bf y}=(y_1,y_2,...,y_{N_h})$ and ${\bf u}=(u_1,u_2,...,u_{N_h})$ to denote their coefficient vectors. Let $U_{ad,h}$ denotes the discretized feasible set, which is defined by
\begin{eqnarray*}
  U_{ad,h}:&=&U_h\cap U_{ad}\\
           %&=&\{z_h\in U_h|~a\leq z_h\leq b, {\rm a.e }\  \mathrm{on}\ \bar{\Omega}_h\} \\
           &=&\left\{z_h=\sum \limits_{i=1}^{N_h} z_i\phi_i(x)~\big{|}~a\leq z_i\leq b, \forall i=1,...,N_h\right\}\subset U_{ad}.
\end{eqnarray*}
Now, a discretized version of the problem (\ref{eqn:orginal problems}) is formulated as follows:
\begin{equation}\label{eqn:discretized problems}
  \left\{ \begin{aligned}
        &\min \limits_{(y_h,u_h)\in Y_h\times U_h}^{}J_h(y_h,u_h)=\frac{1}{2}\|y_h-y_d\|_{L^2(\Omega_h)}^{2}+\frac{\alpha}{2}\|u_h\|_{L^2(\Omega_h)}^{2}+ \beta\|u_h\|_{L^1(\Omega_h)}\\
        &\qquad\quad {\rm{s.t.}}\qquad a(y_h, v_h)=\int_{\Omega}(u_h+y_r)v_h{\rm{d}}x, \qquad  \forall v_h\in Y_h,  \\
          &\qquad \qquad \qquad~  u_h\in U_{ad,h}.
                          \end{aligned} \right.\tag{$\mathrm{P}_{h}$}
 \end{equation}
For the error estimates, we have the following result.
\begin{theorem}{\rm{\textbf{\cite[Proposition 4.3]{WaWa}}}}\label{theorem:error1}
Let us assume that $u^*$ and $u^*_h$ be the optimal control solutions of {\rm(\ref{eqn:orginal problems})} and {\rm(\ref{eqn:discretized problems})}, respectively. Then for every $\alpha_0>0$, $h_0>0$ there exists a constant $C>0$,  such that for all $\alpha\leq\alpha_0$, $h\leq h_0$ the following inequality holds:
\begin{equation}\label{primal-error-estimates}
  \|u^*-u^*_h\|_{L^2(\Omega)}\leq C(\alpha^{-1}h+\alpha^{-3/2}h^2),
\end{equation}
where $C$ is independent of $\alpha, h$.
\end{theorem}

From the perspective of numerical implementation, we introduce the following stiffness and mass matrices
\begin{equation*}
 K_h = \left(a(\phi_i, \phi_j)\right)_{i,j=1}^{N_h},\quad  M_h=\left(\int_{\Omega_h}\phi_i(x)\phi_j(x){\mathrm{d}}x\right)_{i,j=1}^{N_h},
\end{equation*}
and let $y_{r,h}$ and $y_{d,h}$ be the $L^2$-projections of $y_r$ and $y_d$ onto $Y_h$, respectively,
\begin{equation*}
  y_{r,h}=\sum\limits_{i=1}^{N_h}y_r^i\phi_i(x),\quad y_{d,h}=\sum\limits_{i=1}^{N_h}y_d^i\phi_i(x).
\end{equation*}
Similarly, ${\bf y_r}=(y_r^1,y_r^2,...,y_r^{N_h})$ and ${\bf y_d}=(y_d^1,y_d^2,...,y_d^{N_h})$ denote their coefficient vectors, respectively. Then, identifying discretized functions with their coefficient vectors, we can rewrite the problem (\ref{eqn:discretized problems}) in the following way:
\begin{equation}\label{equ:discretized matrix-vector form}
\left\{\begin{aligned}
        &\min\limits_{{\bf y},{\bf u}}^{}~~ \frac{1}{2}\|{\bf y}-{\bf y_{d}}\|_{M_h}^{2}+\frac{\alpha}{2}\|{\bf u}\|_{M_h}^{2}+ \beta\int_{\Omega_h}|\sum_{i=1}^{N_h}{u_i\phi_i(x)}|~\mathrm{d}x\\
        &\ {\rm{s.t.}}\quad K_h{\bf y}=M_h{\bf u}+M_h{\bf y_r},\\
        &\ \quad\quad \ {\bf u}\in [a, b]^{N_h}.%\{v\in\mathbb{R}^n: a\leq v_i\leq b, i=1,...,n\}.
                          \end{aligned} \right.%\tag{${\mathrm{P}}_{h}$}
\end{equation}

\subsection{\textbf {An approximate discretization approach}}
\label{sec:3.2}
To numerically solve problem (\ref{eqn:orginal problems}), a traditional approach is directly solving the primal problem (\ref{equ:discretized matrix-vector form}). However, it is clear that the discretized $L^1$-norm
\begin{equation*}
 \|u_h\|_{L^1(\Omega_h)}= \int_{\Omega_h}|\sum_{i=1}^{N_h}{u_i\phi_i(x)}|\mathrm{d}x
\end{equation*}
is a coupled form with respect to $u_i$ and thus it can not be written as a matrix-vector form. Since its subgradient $\nu_h\in \partial\|u_h\|_{L^1(\Omega_h)}$ will not belong to a finite-dimensional subspace, if directly solving (\ref{equ:discretized matrix-vector form}), it is inevitable to bring some difficulties into numerical calculation. To overcome these difficulties, in \cite{WaWa}, the authors introduced the lumped mass matrix $W_h$
which is a diagonal matrix as:
\begin{equation}\label{def-lumped-mass}
  W_h:={\rm{diag}}\left(\int_{\Omega_h}\phi_i(x)\mathrm{d}x\right)_{i=1}^{N_h},
\end{equation}
and defined an alternative discretization of the $L^1$-norm:
\begin{equation}\label{equ:approximal L1}
  \|u_h\|_{L^{1}_h(\Omega)}:=\sum_{i=1}^{N_h}|u_i|\int_{\Omega_h}\phi_i(x)\mathrm{d}x=\|W_h{\bf u}\|_1,
\end{equation}
which is a weighted $l^1$-norm of the coefficients of $u_h$.
More importantly, the following results about the mass matrix $M_h$ and the lumped mass matrix $W_h$ hold.
\begin{proposition}\label{eqn:martix properties1}
$\forall$ ${\bf z}=(z_1,z_2,...,z_{N_h})\in \mathbb{R}^{N_h}$, the following inequalities hold:
\begin{eqnarray}
 \label{Winequality1}&&\|{\bf z}\|^2_{M_h}\leq\|{\bf z}\|^2_{W_h}\leq \gamma\|{\bf z}\|^2_{M_h} \quad where \ \gamma=
 \left\{ \begin{aligned}
         &4  \quad if \ n=2, \\
         &5  \quad if \ n=3,
                           \end{aligned} \right.
                           \\
 \label{Winequality2}&&\int_{\Omega_h}|\sum_{i=1}^n{z_i\phi_i(x)}|~\mathrm{d}x\leq\|W_h{\bf z}\|_1.
\end{eqnarray}
\proof
Based on non-negativity and partition of unity of the nodal basis functions, utilizing convexity argument, it is easy to obtain
\begin{equation*}
 \left(\sum\limits_{i=1}^{n}z_i\phi_i(x)\right)^2\leq\sum\limits_{i=1}^{n}(z_i)^2\phi_i(x),\quad
 |\sum\limits_{i=1}^{n}z_i\phi_i(x)|\leq\sum\limits_{i=1}^{n}|z_i|\phi_i(x).
\end{equation*}
This implies
\begin{equation*}
\begin{aligned}
&\|{\bf z}\|^2_{M_h}=\int_{\Omega_h}\left(\sum\limits_{i=1}^{n}z_i\phi_i(x)\right)^2{\rm d}x\leq\int_{\Omega_h}\sum\limits_{i=1}^{n}(z_i)^2\phi_i(x){\rm d}x=\|{\bf z}\|^2_{W_h},\\
&\int_{\Omega_h}|\sum\limits_{i=1}^{n}z_i\phi_i(x)|{\rm d}x\leq\int_{\Omega_h}\sum\limits_{i=1}^{n}|z_i|\phi_i(x){\rm d}x=\|W_h{\bf z}\|_1.
\end{aligned}
\end{equation*}
For a proof of the inequality $\|{\bf z}\|^2_{W_h}\leq \gamma\|{\bf z}\|^2_{M_h}$, we refer to {\rm\cite[Table 1]{Wathen}}.
\endproof

\end{proposition}

Thus, we provide a discretization of problem {\rm(\ref{eqn:orginal problems})}:
\begin{equation}\label{eqn:approx discretized problems}
  \left\{ \begin{aligned}
        &\min \limits_{(y_h,u_h)\in Y_h\times U_h}^{}J_h(y_h,u_h)=\frac{1}{2}\|y_h-y_d\|_{L^2(\Omega_h)}^{2}+\frac{\alpha}{2}\|u_h\|_{L^2(\Omega_h)}^{2}+ \beta\|u_h\|_{L_h^1(\Omega_h)}\\
        &\qquad\quad {\rm{s.t.}}\qquad a(y_h, v_h)=\int_{\Omega}(u_h+y_r)v_h{\rm{d}}x, \qquad  \forall v_h\in Y_h,  \\
          &\qquad \qquad \qquad~  u_h\in U_{ad,h}.
                          \end{aligned} \right.\tag{$\widetilde{\mathrm{P}}_{h}$}
 \end{equation}
where $\|\cdot\|_{L_h^1(\Omega_h)}$ is defined in (\ref{equ:approximal L1}). Similarly, (\ref{eqn:approx discretized problems}) can also be rewritten as the following matrix-vector form:
\begin{equation}\label{equ:approx discretized matrix-vector form}
\left\{\begin{aligned}
        &\min\limits_{{\bf y,u}}^{}~~\frac{1}{2}\|{\bf y}-{\bf y_{d}}\|_{M_h}^{2}+\frac{\alpha}{2}\|{\bf u}\|_{M_h}^{2}+ \beta\|W_h {\bf u}\|_{1}\\
        &\ {\rm{s.t.}}\quad K_h{\bf y}=M_h ({\bf u} +{\bf y_r}),\\
        &\ \quad\quad \ {\bf u}\in [a,b]^{N_h}.%\{v\in\mathbb{R}^n: a\leq v_i\leq b, i=1,...,n\}.
                          \end{aligned} \right.%\tag{$\widetilde{\mathrm{P}}_{h}$}
\end{equation}

About the error estimates results between (\ref{eqn:orginal problems}) and (\ref{eqn:approx discretized problems}), we have the following result, see \cite[Corollary 4.6]{WaWa} for more details.
\begin{theorem}{\rm{\textbf{\cite[Corollary 4.6]{WaWa}}}}\label{theorem:error1}
Let us assume that $u^*$ and $u^*_h$ be the optimal control solutions of {\rm(\ref{eqn:orginal problems})} and {\rm(\ref{eqn:approx discretized problems})}, respectively. Then for every $\alpha_0>0$, $h_0>0$ there exists a constant $C>0$,  such that for all $\alpha\leq\alpha_0$, $h\leq h_0$ the following inequality holds
\begin{equation*}
  \|u^*-u^*_h\|_{L^2(\Omega)}\leq C(\alpha^{-1}h+\alpha^{-3/2}h^2),
\end{equation*}
where $C$ is independent of $\alpha, h$.
\end{theorem}

As it turned out, the additional approximation step (\ref{equ:approximal L1}) would not disturb the convergence estimate, in fact, both the error orders $h$ and $\alpha$ in the estimate remain unchanged. However, the approximation of $L^1$-norm (\ref{equ:approximal L1}) inevitably brings additional error. Thus, it is necessary to analyze the error between $\|u_h\|_{L^{1}_h(\Omega)}$ and $\|u_h\|_{L^{1}(\Omega)}$.

To achieve our goal, let us first introduce the nodal interpolation operator $I_h$. For a given regular and quasi-uniform triangulation $\mathcal{T}_h$ of $\Omega$ with nodes $\{x_i\}_{i=1}^{N_h}$, we define
\begin{equation}\label{nodal interpolation operator}
(I_hw)(x)=\sum_{i=1}^{N_h}w(x_i)\phi_i(x) \ {\rm\ for\ any}\ w\in L^1(\Omega).
\end{equation}
About the interpolation error estimate, we have the following result, see  {\rm\cite[Theorem 3.1.6]{Ciarlet}} for more details.

\begin{lemma}\label{interpolation error estimate}
For all $w\in W^{k+1,p}(\Omega)$, $k\geq 0$, $p,q\in [0,+\infty)$, and $0\leq m\leq k+1$, we have
\begin{equation}\label{interpolation error estimate inequilaty}
  \|w-I_hw\|_{W^{m,q}(\Omega)}\leq c_I h^{k+1-m}\|w\|_{W^{k+1,p}(\Omega)}.
\end{equation}
\end{lemma}
Thus, according to Lemma \ref{interpolation error estimate}, we have the following error estimate results.
\begin{proposition}\label{eqn:martix properties}
$\forall$ ${\bf z}=(z_1,z_2,...,z_{N_h})\in \mathbb{R}^{N_h}$, let $z_h=\sum\limits_{i=1}^{N_h}z_i\phi_i$, then the following inequalities hold
\begin{eqnarray}
 \label{Winequality3}&&0\leq\|{\bf z}\|^2_{W_h}-\|{\bf z}\|^2_{M_h}\leq C\|z_h^2\|_{H^2(\Omega)}h^2,\\%\mathcal{O}(h^2)\\
 \label{Winequality4}&&0\leq\|W_h{\bf z}\|_1-\int_{\Omega_h}|\sum_{i=1}^n{z_i\phi_i(x)}|~\mathrm{d}x\leq C\|z_h\|_{H^1(\Omega)}h,
\end{eqnarray}
where $C$ is a constant.
\proof
First, we have
\begin{equation*}
  \|{\bf z}\|^2_{W_h}=\int_{\Omega_h}\sum_{i=1}^{N_h}(z_i)^2\phi_i(x)\mathrm{d}x=\int_{\Omega_h}I_h(z_h)^2(x)\mathrm{d}x,
\end{equation*}
where $I_h$ is the nodal interpolation operator. Since $z_h\in \mathcal{P}_1\subset W^{1,2}(\Omega)$, we have $(z_h)^2\in W^{2,2}(\Omega)$. Thus by Lemma {\rm\ref{interpolation error estimate}}, we get

\begin{eqnarray*}
  \|{\bf z}\|^2_{W_h}-\|{\bf z}\|^2_{M_h}&=&\int_{\Omega_h}I_h(z_h)^2(x)\mathrm{d}x-\int_{\Omega_h}(z_h)^2(x)\mathrm{d}x\\
  &=&\|I_h(z_h)^2- (z_h)^2\|_{L^1(\Omega)}\\
  &\leq& c_{\Omega}\|I_h(z_h)^2- (z_h)^2\|_{L^2(\Omega)}\\
  &\leq&c_{\Omega}c_I\|z_h^2\|_{H^2(\Omega)}h^2\\
  &=&C\|z_h^2\|_{H^2(\Omega)}h^2.
\end{eqnarray*}
Similarly, we have
\begin{eqnarray*}
  \|W_h{\bf z}\|_1-\int_{\Omega_h}|\sum_{i=1}^n{z_i\phi_i(x)}|~\mathrm{d}x
  &=& \|I_h|z_h|\|_{L^1(\Omega)}- \|z_h\|_{L^1(\Omega)}\\
  &=&\int_{\Omega_h}I_h|z_h|(x)-|z_h(x)|\mathrm{d}x\\
  &\leq& c_{\Omega}\|I_h|z_h|- |z_h|\|_{L^2(\Omega)}\\
  &\leq&c_{\Omega}c_I\|z_h\|_{H^1(\Omega)}h\\
  &=&C\|z_h\|_{H^1(\Omega)}h\\
\end{eqnarray*}
where the last equation is due to $|z_h|\in W^{1,2}(\Omega)$.
\endproof
\end{proposition}

\section{Duality-based approach}
\label{sec:4}
As we said, in \cite{mABCDSOPT}, the authors considered using the duality-based approach to solve (\ref{eqn:orginal problems}). Thus, in this section, we will introduce the dual problem of (\ref{eqn:orginal problems}) and give a brief sketch of the symmetric Gauss-Seidel based majorized ABCD (sGS-mABCD) method for the dual problem. At last, in order to achieve our ultimate goal of finding the optimal control and the optimal state, we will introduce the primal problem of the discretized dual problem.

\subsection{\textbf{Dual problem of (\ref{eqn:orginal problems})}}
About the dual problem of (\ref{eqn:orginal problems}), we have the following proposition.
\begin{proposition}\label{primal problem of Dh}
The dual problem of {\rm(\ref{eqn:orginal problems})} can be written, in its equivalent minimization form, as
\begin{equation}\label{eqn:dual problem}
\begin{aligned}
\min\ \Phi(\lambda,p,\mu):=&{\frac{1}{2}\|A^*p-y_d\|_{L^2(\Omega)}^2}+ \frac{1}{2\alpha}\|-p+\lambda+\mu\|_{L^2(\Omega)}^2+\langle p, y_r\rangle_{L^2(\Omega)}\\
&+\delta_{\beta B_{\infty}(0)}(\lambda)+\delta^*_{U_{ad}}(\mu)-\frac{1}{2}\|y_d\|_{L^2(\Omega)}^2,\end{aligned}\tag{$\mathrm{D}$}
\end{equation}
where $p\in H^1_0(\Omega)$, $\lambda,\mu\in L^2(\Omega)$, $B_{\infty}(0):=\{\lambda\in L^2(\Omega): \|\lambda\|_{L^\infty(\Omega)}\leq 1\}$, and for any given nonempty, closed convex subset $C$ of $L^2(\Omega)$, $\delta_{C}(\cdot)$ is the indicator function of $C$. Based on the $L^2$-inner product, the conjugate of $\delta_{C}(\cdot)$ is defined as follows
\begin{equation*}
 \delta^*_{C}(w^*)=\sup\limits_{w\in C}^{}{\langle w^*,w\rangle}_{L^2(\Omega)}.
\end{equation*}
\proof
{Firstly, by introducing two artificial variables $v\in{L^2(\Omega)}$ and $w\in{L^2(\Omega)}$, we can rewrite {\rm(\ref{eqn:orginal problems})} as:
\begin{equation}\label{equ:decoupled new form}
\left\{\begin{aligned}
        &\min\limits_{y,u,v,w}^{}~~\bar{J}(y,u,v,w)= \frac{1}{2}\|y-y_d\|_{L^2(\Omega)}^{2}+\frac{\alpha}{2}\|u\|_{L^2(\Omega)}^{2}+\beta\|v\|_{L^1(\Omega)} +\delta_{[a,b]}(w)\\
        & {\rm{s.t.}}\qquad\quad Ay=B (u +y_r),\\
        &\ \qquad \qquad u-v=0,\\
        &\ \qquad \qquad u-w=0.
                          \end{aligned} \right.%\tag{$\widetilde{\mathrm{P}}_{h}$}
\end{equation}
Considering the Lagrangian function associated with {\rm(\ref{equ:decoupled new form})}, we have
\begin{eqnarray}\label{Lagrangian function for dual}
\begin{aligned}
  L(y,u,v,w;p,\lambda,\mu)&=\bar{J}(y,u,v,w)+\langle p,Ay-B (u +y_r)\rangle_{L^2(\Omega)}+\langle\lambda,u-v\rangle_{L^2(\Omega)}+\langle\mu,u-w\rangle_{L^2(\Omega)}.
\end{aligned}
\end{eqnarray}
Now, we can derive
\begin{eqnarray}\label{equ:minimum Lagrangian function}
\left\{\begin{aligned}
  \inf_y L(y,u,v,w;p,\lambda,\mu) &= \inf_y{\frac{1}{2}\|y-y_{d}\|_{L^2(\Omega)}^{2}}+\langle p,Ay\rangle =-\frac{1}{2}\|A^* p-y_{d}\|_{L^2(\Omega)}^2+\frac{1}{2}\|y_d\|^2_{L^2(\Omega)},\\
  \inf_u L(y,u,v,w;p,\lambda,\mu)&=\inf_u \frac{\alpha}{2}\|u\|^2_{L^2(\Omega)}-\langle p,B u\rangle+\langle\lambda,u\rangle_{L^2(\Omega)}+\langle\mu,u\rangle_{L^2(\Omega)}\\
  &=-\frac{1}{2\alpha}\|\lambda+ \mu- p\|_{L^2(\Omega)}^2,\\
   \inf_vL(y,u,v,w;p,\lambda,\mu)&=\inf_v-\langle\lambda,v\rangle_{L^2(\Omega)}+\beta\|v\|_1=-\delta_{[-\beta,\beta]}(\lambda), \\
  \inf_wL(y,u,v,w;p,\lambda,\mu)&=\inf_w-\langle\mu,w\rangle_{L^2(\Omega)}+\delta_{[a,b]}(w)=- \delta^*_{[a,b]}(\mu).
\end{aligned}\right.
\end{eqnarray}
Thus,
\begin{eqnarray*}
\begin{aligned}
\min_{y,u,v,w}L(y,u,v,w;p,\lambda,\mu)=&-\frac{1}{2}\|A^* p-y_{d}\|_{L^2(\Omega)}^2- \frac{1}{2\alpha}\|\lambda+ \mu- p\|_{L^2(\Omega)}^2-\langle y_r, p\rangle_{L^2(\Omega)}\\
  &-\delta_{[-\beta,\beta]}(\lambda)-\delta^*_{[a,b]}({}\mu)+\frac{1}{2}\|y_d\|^2_{L^2(\Omega)},
\end{aligned}
\end{eqnarray*}
and $\max_{p,\lambda,\mu}\min_{y,u,v,w}L(y,u,v,w;p,\lambda,\mu)$ is an equivalent maximization form of the dual problem {\rm(\ref{eqn:dual problem})}. Thus, we complete the proof.}
\endproof
\end{proposition}

Employing the piecewise linear, globally continuous finite elements to discretize all the dual variables, then
a type of finite element discretization of (\ref{eqn:dual problem}) is given as follows
\begin{equation}\label{eqn:discretized matrix-vector dual problem}
\begin{aligned}
\min\limits_{{\bm \lambda},{\bf p},{\bm \mu}\in \mathbb{R}^{N_h}}
\Phi_h({\bm \lambda},{\bf p},{\bm \mu}):=&
\frac{1}{2}\|K_h {\bf p}-{M_h}{\bf y_{d}}\|_{M_h^{-1}}^2+ \frac{1}{2\alpha}\|{\bm \lambda}+{\bm \mu}-{\bf p}\|_{M_h}^2+\langle M_h{\bf y_r}, {\bf p}\rangle\\
&+ \delta_{[-\beta,\beta]}({\bm\lambda})+ \delta^*_{[a,b]}({M_h}{\bm\mu})-\frac{1}{2}\|{\bf y_d}\|^2_{M_h}.
\end{aligned}\tag{$\mathrm{D_h}$}
\end{equation}

\subsection{\textbf{A sGS based majorized ABCD method for (\ref{eqn:discretized matrix-vector dual problem})}}\label{sec:4.2}
Obviously, (\ref{eqn:discretized matrix-vector dual problem}) belongs to a general class of unconstrained, multi-block convex optimization problems with coupled objective function, that is
\begin{equation}\label{eqn:model problem}
\begin{aligned}
\min_{v, w} \theta(v,w):= f(v)+ g(w)+ \phi(v, w),
\end{aligned}
\end{equation}
where $f: \mathcal{V}\rightarrow (-\infty, +\infty ]$ and $g: \mathcal{W}\rightarrow  (-\infty, +\infty ]$ are two convex functions (possibly nonsmooth), $\phi: \mathcal{V}\times \mathcal{W}\rightarrow  (-\infty, +\infty ]$ is a smooth convex function, and $\mathcal{V}$, $\mathcal{W}$ are real finite dimensional Hilbert spaces.

Taking advantage of the structure of (\ref{eqn:model problem}), in \cite[Chapter 3]{CuiYing}, Cui proposed an inexact majorized accelerated block coordinate descent (imABCD) method for solving it.
Under suitable assumptions and certain inexactness criteria, the author can prove that the imABCD method enjoys the impressive $O(1/k^2)$ iteration complexity. In \cite{mABCDSOPT}, which is inspired by the success of the imABCD method, the authors combine the virtues of the recent advances in the inexact sGS technique and the imABCD method and propose a sGS based majorized ABCD method (called sGS-mABCD) to efficiently and fast solve problem (\ref{eqn:discretized matrix-vector dual problem}).

In this paper, we give a brief sketch of the sGS-mABCD method. First, we express (\ref{eqn:discretized matrix-vector dual problem}) in the form of (\ref{eqn:model problem}) with
%$v=\bm\mu$, $w=(\bm\lambda, {\bf p})$ and
%\begin{eqnarray}
%% \nonumber to remove numbering (before each equation)
%  f(v) &=&\delta^*_{[a,b]}({M_h}\bm\mu)\label{f function for Dh},\\
%  g(w) &=& \delta_{[-\beta,\beta]}(\bm\lambda)+\frac{1}{2}\|K_h {\bf p}-{M_h}{\bf y_{d}}\|_{M_h^{-1}}^2+\langle M_h {\bf y_r}, {\bf p}\rangle-\frac{1}{2}\|{\bf y_d}\|^2_{M_h}\label{g function for Dh}, \\
%  \phi(v, w) &=& \frac{1}{2\alpha}\|\bm\lambda+ \bm\mu- {\bf p}\|_{M_h}^2\label{phi function for Dh}.
%\end{eqnarray}
%In the second variant, we exchange the iteration order of $\bm \mu$ and $(\bm \lambda, {\bf p})$, which implies we choose
$v=(\bm \lambda, {\bf p})$, $w=\bm \mu$ and
\begin{eqnarray}
  f(v) &=&\delta_{[-\beta,\beta]}(\bm\lambda)+\frac{1}{2}\|K_h {\bf p}-{M_h}{\bf y_{d}}\|_{M_h^{-1}}^2+\langle M_h {\bf y_r}, {\bf p}\rangle-\frac{1}{2}\|{\bf y_d}\|^2_{M_h} \label{g function for Dh}, \\
  g(w) &=&\delta^*_{[a,b]}({M_h}\bm\mu)\label{f function for Dh},\\
  \phi(v, w) &=&\frac{1}{2\alpha}\|\bm\lambda+ \bm\mu- {\bf p}\|_{M_h}^2\label{phi function for Dh}.
\end{eqnarray}

The detailed framework of the sGS-mABCD method for (\ref{eqn:discretized matrix-vector dual problem}) is given as follows. It should be stressed that a block symmetric Gauss-Seidel decomposition technique for
convex composite quadratic programming plays a key role in solving the $(\bm \lambda, {\bf p})$-subproblem.

\begin{algorithm}[H]\label{algo1:Full inexact ABCD algorithm for (Dh)}
  \caption{\textbf{A sGS-mABCD method for (\ref{eqn:discretized matrix-vector dual problem})}}
\begin{description}
\item[]
  \item [\bf Input]$(\tilde{\bm \lambda}^1, \tilde{{\bf p}}^1,\tilde{\bm \mu}^1)=({\bm \lambda}^0, {{\bf p}}^0,\bm \mu^0)\in  [-\beta,\beta]\times \mathbb{R}^{N_h}\times{\rm dom} (\delta^*_{[a,b]})$. Set $k= 1, t_1= 1.$
  \item [\bf Output]$ ({\bm \lambda}^k, {{\bf p}}^k,{\bm \mu}^k)$
  \item [\bf Iterate until convergence]
  \end{description}
\begin{description}
  \item[\bf Step 1] Utilizing the block symmetric Gauss-Seidel iteration to compute block-$(\bm \lambda^k,{\bf p}^k)$ as follows:
  \begin{itemize}
       \item[$\bullet$] (Backward GS sweep) Compute
       \begin{equation*}
       \hat{{\bf p}}^{k}=\arg\min\frac{1}{2}\|K_h {\bf p}-{M_h}{\bf y_{d}}\|_{M_h^{-1}}^2+ \frac{1}{2\alpha}\|{\bf p}-\tilde{\bm \lambda}^k-\tilde{\bm \mu}^k\|_{M_h}^2+\langle M_h {\bf y_r}, {\bf p}\rangle,
       \end{equation*}
       \item[$\bullet$] Then compute
       \begin{equation*}
         {\bm \lambda}^{k}
        =\arg\min\delta_{[-\beta,\beta]}(\bm \lambda)+\frac{1}{2\alpha}\|\bm \lambda-(\hat{{\bf p}}^{k}-\tilde{\bm \mu}^k)\|_{M_h}^2+\frac{1}{2\alpha}\|\bm \lambda- \tilde{\bm \lambda}^{k}\|_{W_h-M_h}^2,
       \end{equation*}
       \item[$\bullet$] (Forward GS sweep) Compute
       \begin{equation*}
         {{\bf p}}^{k}=\arg\min\frac{1}{2}\|K_h {\bf p}-{M_h}{\bf y_{d}}\|_{M_h^{-1}}^2+ \frac{1}{2\alpha}\|{\bf p}-{\bm \lambda}^k-\tilde{\bm \mu}^k\|_{M_h}^2+\langle M_h {\bf y_r}, {\bf p}\rangle,
       \end{equation*}
  \end{itemize}
 \item[\bf Step 2] Compute block-$\bm \mu$
\begin{equation*}
{\bm \mu}^{k}=\arg\min\delta^*_{[a,b]}(M_h\bm \mu)+\frac{1}{2\alpha}\|\bm \mu-({\bf p}^k-\bm \lambda^k)\|_{M_h}^2+\frac{1}{2\alpha}\|\bm \mu-\tilde{\bm \mu}^k\|^2_{\gamma M_hW_h^{-1}M_h-M_h},
\end{equation*}
where $\gamma$ is defined in Proposition {\rm\ref{eqn:martix properties1}}.
  \item[\bf Step 3] Set $t_{k+1}=\frac{1+\sqrt{1+4t_k^2}}{2}$ and $\beta_k=\frac{t_k-1}{t_{k+1}}$, Compute
\begin{eqnarray*}
\tilde{\bm \lambda}^{k+1}= {\bm \lambda}^{k}+ \beta_{k}({\bm \lambda}^{k}-{\bm \lambda}^{k-1}),\quad
\tilde{{\bf p}}^{k+1}={{\bf p}}^{k}+ \beta_{k}({{\bf p}}^{k}-{{\bf p}}^{k-1}), \quad\tilde{\bm \mu}^{k+1}={\bm \mu}^{k}+ \beta_{k}({\bm \mu}^{k}-{\bm \mu}^{k-1}).
\end{eqnarray*}
\end{description}
\end{algorithm}

About the iteration complexity of Algorithm \ref{algo1:Full inexact ABCD algorithm for (Dh)}, we have the following results. For more details, one can refer to \cite[Theorem 7]{meshindependencemABCD}.
%of the objective values.
\begin{proposition}\label{sGS-imABCD convergence}
Suppose that the solution set $\Omega$ of the problem {\rm(\ref{eqn:discretized matrix-vector dual problem})} is non-empty. Let
${\bf z}^*=(\bm \lambda^*,{\bf p}^*,\bm \mu^*)\in \Omega$. Let $\{{\bf z}^k\}:=\{({\bm \lambda}^k,{{\bf p}}^k,{\bm \mu}^k)\}$ be the sequence generated by the Algorithm {\rm\ref{algo1:Full inexact ABCD algorithm for (Dh)}}. Then we have
\begin{equation}\label{iteration complexity of discretized dual problem2}
\Phi_h({\bf z}^k)- \Phi_h({\bf z}^*)\leq \frac{4\tau^1_h}{(k+1)^2} \quad \forall k\geq 1,
\end{equation}
where $\Phi(\cdot)$ is the objective function of the dual problem {\rm(\ref{eqn:discretized matrix-vector dual problem})} and
\begin{eqnarray}
  &&\tau^1_h=\frac{1}{2}\|{\bf z}^0- {\bf z}^*\|_{\mathcal{S}^1_h}^1,\label{convergence factor21} \\
  &&\mathcal{S}^1_h=\frac{1}{\alpha}\left(
                                      \begin{array}{ccc}
                                        M_hG_h^{-1}M_h+W_h-M_h &  \ 0&\  0\\
                                        0 &  \ 0&  \ 0\\
                                        0&  \ 0&  \ \gamma M_hW_h^{-1}M_h\\
                                      \end{array}
                                    \right)\succeq0\label{convergence factor22}.\\
  &&G_h=M_h+\alpha K_h M_h^{-1}K_h\label{convergence factor23}.
\end{eqnarray}
\end{proposition}

Next, another key issue should be considered is how measures of the convergence behavior of the iteration sequence vary with the level of approximation. In other words, we should analyse whether the ``discretized" convergence factor $\tau^1_h$ could be uniformly bounded by a constant which is independent of the mesh size $h$. In order to show these results, let us first present some bounds on the Rayleigh quotients of $K_h$ and $M_h$, one can see Proposition 1.29 and Theorem 1.32 in \cite{spectralproperty} for more details.

\begin{lemma}\label{spectral property}
For $\mathcal{P}1$ approximation on a regular and quasi-uniform subdivision of $\mathbb{R}^n$ which satisfies Assumption {\rm\ref{regular and quasi-uniform triangulations}}, and for any ${\bf x}\in \mathbb{R}^{N_h}$, the mass matrix $M_h$ approximates the scaled identity matrix in the sense that
\begin{equation*}
c_1 h^2\leq \frac{{\bf x}^{T}M_h{\bf x}}{{\bf x}^{T}{\bf x}}\leq c_2 h^2 \ if\  n=2, \ {\rm and}\ c_1 h^3\leq \frac{{\bf x}^{T}M_h{\bf x}}{{\bf x}^{T}{\bf x}}\leq c_2 h^3 \ if\  n=3,
\end{equation*}
the stiffness matrix $K_h$ satisfies
\begin{equation*}
d_1h^2\leq \frac{{\bf x}^{T}K_h{\bf x}}{{\bf x}^{T}{\bf x}}\leq d_2 \ if\  n=2, \ {\rm and }\ d_1h^3\leq \frac{{\bf x}^{T}K_h{\bf x}}{{\bf x}^{T}{\bf x}}\leq d_2 h \ if\  n=3,
\end{equation*}
where the constants $c_1$, $c_2$, $d_1$ and $d_2$ are independent of the mesh size $h$.
\end{lemma}

Based on Lemma \ref{spectral property}, we can easily obtain the following lemma.

\begin{lemma}\label{spectral property of Gh}
For any ${\bf x}\in \mathbb{R}^{N_h}$, there exist four constants $u_1$, $u_2$, $l_1$, $l_2$ and $h_0>0$, such that for any $0<h<h_0$, the matrix $G_h$ which defined in {\rm(\ref{convergence factor23})} satisfies the following inequalities
\begin{equation}\label{spectral property of Gh12}
  \begin{aligned}
&l_1 h^2\leq \frac{{\bf x}^{T}G_h{\bf x}}{{\bf x}^{T}{\bf x}}\leq u_1 \frac{1}{h^2} \quad if\  n=2, \quad
l_2 h^3\leq \frac{{\bf x}^{T}G_h{\bf x}}{{\bf x}^{T}{\bf x}}\leq u_2 \frac{1}{h} \quad\ if\  n=3.
  \end{aligned}
\end{equation}
%\begin{eqnarray}
%&&l_1 h^2\leq \frac{{\bf x}^{T}G_h{\bf x}}{{\bf x}^{T}{\bf x}}\leq u_1 \frac{1}{h^2} \quad if\  n=2 \label{spectral property of Gh1},\\
%&&l_2 h^3\leq \frac{{\bf x}^{T}G_h{\bf x}}{{\bf x}^{T}{\bf x}}\leq u_2 \frac{1}{h} \quad\ if\  n=3\label{spectral property of Gh2}.
%\end{eqnarray}
\end{lemma}

At first glance, it appears that the largest eigenvalue of the matrix $\mathcal{S}^1_h$ can not be uniformly bounded by a constant $C$ for the reason of $G_h$. However, based on Lemma \ref{spectral property} and Lemma \ref{spectral property of Gh}, it is easy to prove that there exists $h_0>0$, such that for any $0<h<h_0$, the matrix $M_hG_h^{-1}M_h$ satisfies the following properties
\begin{equation*}\label{property of MinvGM}
\begin{aligned}
 &\lambda_{\max}(M_hG_h^{-1}M_h)=O(h^2) \quad {\rm for}\ n=2,\quad \lambda_{\max}(M_hG_h^{-1}M_h)=O(h^3) \quad {\rm for}\ n=3,
\end{aligned}
\end{equation*}
where $ \lambda_{\max}(\cdot)$ represents the largest eigenvalue of a given matrix. Furthermore, we have
\begin{equation}\label{property of Sh}
\begin{aligned}
  \lambda_{\max}(\mathcal{S}^1_h)&=\frac{1}{\alpha}\max\{ \lambda_{\max}(M_hG_h^{-1}M_h+W_h-M_h),\lambda_{\max}(\gamma M_hW_h^{-1}M_h)\}\\
  &=\left\{\begin{aligned}
      & O(h^2) \quad {\rm for}\ n=2,\\
      & O(h^3) \quad {\rm for}\ n=3.
    \end{aligned}\right.
\end{aligned}.
\end{equation}
In other words, we can say that the largest eigenvalue of the matrix $\mathcal{S}^1_h$ can be uniformly bounded by a constant $C$ which is independent of the mesh size $h$, which implies the ``discretized" convergence factor $\tau^1_h$ could be uniformly bounded by a constant. From this point of view, the mesh independence of Algorithm sGS-mABCD is its another advantage.

\subsection{\textbf{A majorized ABCD method with semismooth Newton for (\ref{eqn:discretized matrix-vector dual problem})}}\label{sec:4.3}
If we carefully check Algorithm \ref{algo1:Full inexact ABCD algorithm for (Dh)}, it should be pointed out that the information of the accelerated points $\{\tilde {{\bf p}}^k\}_{k>1}$ is not used in the whole iterative process. Although in theory such iterative scheme not affect the convergence result as shown in Proposition \ref{sGS-imABCD convergence}, the lack of such acceleration information may affect the actual convergence rate of the algorithm in the numerical implementation. In order to more efficiently achieve a high accuracy, we give a majorized ABCD method with semismooth Newton conjugate gradient algorithm for (\ref{eqn:discretized matrix-vector dual problem}). Specifically, instead of using the symmetric Gauss-Seidel technique to solve the $({\bm \lambda}, {\bf p})$-subproblem, we employ a semismooth Newton conjugate gradient (SNCG) algorithm introduced in \cite{SNCG1,SNCG2} to solve it.

The detailed framework of the majorized ABCD with semismooth Newton conjugate gradient algorithm is presented as follows.

\begin{algorithm}[H]\label{algo2:Full inexact ABCD algorithm for (Dh)}
  \caption{\textbf{An majorized ABCD-SNCG for (\ref{eqn:discretized matrix-vector dual problem})}}
\begin{description}
\item[]
  \item [\bf Input]$(\tilde{\bm \lambda}^1, \tilde{{\bf p}}^1,\tilde{\bm \mu}^1)=({\bm \lambda}^0, {{\bf p}}^0,\bm \mu^0)\in  [-\beta,\beta]\times \mathbb{R}^{N_h}\times{\rm dom} (\delta^*_{[a,b]})$. Set $k= 1, t_1= 1.$
  \item [\bf Output]$ ({\bm \lambda}^k, {{\bf p}}^k,{\bm \mu}^k)$
  \item [\bf Iterate until convergence]
  \end{description}
\begin{description}
  \item[\bf Step 1] Utilizing the semismooth Newton-CG algorithm to compute block-$(\bm \lambda^k,{\bf p}^k)$ as follows:
       \begin{equation*}
       \begin{aligned}
         ({\bm \lambda^k},{\bf p}^{k})=\arg\min&\frac{1}{2\alpha}\|{\bf p}-{\bm \lambda}-\tilde{\bm \mu}^k\|_{M_h}^2+\frac{1}{2}\|K_h {\bf p}-{M_h}{\bf y_{d}}\|_{M_h^{-1}}^2+\langle M_h {\bf y_r}, {\bf p}\rangle\\
         &+\delta_{[-\beta,\beta]}(\bm \lambda)+\frac{1}{2\alpha}\|\bm \lambda- \tilde{\bm \lambda}^{k}\|_{W_h-M_h}^2,
       \end{aligned}
       \end{equation*}
  which is equivalent to using the semismooth Newton-CG algorithm to solve
\begin{equation*}
         \begin{aligned}
         {\bf p}^{k}=\arg\min& \frac{1}{2\alpha}\|{\bf p}-{\bm \lambda({\bf p})}-\tilde{\bm \mu}^k\|_{M_h}^2+\frac{1}{2}\|K_h {\bf p}-{M_h}{\bf y_{d}}\|_{M_h^{-1}}^2+ \langle M_h {\bf y_r}, {\bf p}\rangle\\
=\arg\min& \frac{1}{2\alpha}\|{\bf p}-{\rm\Pi}_{[-\beta,\beta]}(\tilde{\bm \lambda}^{k}+W_h^{-1}M_h({\bf p}-\tilde{\bm \mu}^k-\tilde{\bm \lambda}^{k}))-\tilde{\bm \mu}^k\|_{M_h}^2\\
&+\frac{1}{2}\|K_h {\bf p}-{M_h}{\bf y_{d}}\|_{M_h^{-1}}^2+\langle M_h {\bf y_r}, {\bf p}\rangle.
         \end{aligned}
       \end{equation*}
Then
\begin{equation*}
         \begin{aligned}
  &{\bm \lambda}^{k}={\rm\Pi}_{[-\beta,\beta]}(\tilde{\bm \lambda}^{k}+W_h^{-1}M_h({\bf p}^k-\tilde{\bm \mu}^k-\tilde{\bm \lambda}^{k})).
         \end{aligned}
       \end{equation*}
 \item[\bf Step 2] Compute block-$\bm \mu$
\begin{equation*}
{\bm \mu}^{k}=\arg\min\delta^*_{[a,b]}(M_h\bm \mu)+\frac{1}{2\alpha}\|\bm \mu-({\bf p}^k-\bm \lambda^k)\|_{M_h}^2+\frac{1}{2\alpha}\|\bm \mu-\tilde{\bm \mu}^k\|^2_{\gamma M_hW_h^{-1}M_h-M_h},
\end{equation*}
where $\gamma$ is defined in Proposition {\rm\ref{eqn:martix properties1}}.
  \item[\bf Step 3] Set $t_{k+1}=\frac{1+\sqrt{1+4t_k^2}}{2}$ and $\beta_k=\frac{t_k-1}{t_{k+1}}$, Compute
\begin{eqnarray*}
\tilde{\bm \lambda}^{k+1}= {\bm \lambda}^{k}+ \beta_{k}({\bm \lambda}^{k}-{\bm \lambda}^{k-1}),\quad
\tilde{{\bf p}}^{k+1}={{\bf p}}^{k}+ \beta_{k}({{\bf p}}^{k}-{{\bf p}}^{k-1}), \quad\tilde{\bm \mu}^{k+1}={\bm \mu}^{k}+ \beta_{k}({\bm \mu}^{k}-{\bm \mu}^{k-1}).
\end{eqnarray*}
\end{description}
\end{algorithm}

Similarly, about the iteration complexity of Algorithm \ref{algo2:Full inexact ABCD algorithm for (Dh)}, we have the following result.
\begin{proposition}\label{sGS-imABCD convergence}
Suppose that the solution set $\Omega$ of the problem {\rm(\ref{eqn:discretized matrix-vector dual problem})} is non-empty. Let
${\bf z}^*=(\bm \lambda^*,{\bf p}^*,\bm \mu^*)\in \Omega$. Let $\{{\bf z}^k\}:=\{({\bm \lambda}^k,{\bf p}^k),{\bm \mu}^k\}$ be the sequence generated by the Algorithm {\rm\ref{algo2:Full inexact ABCD algorithm for (Dh)}}. Then we have
\begin{equation}\label{iteration complexity of discretized dual problem}
\Phi_h({\bf z}^k)- \Phi_h({\bf z}^*)\leq \frac{4\tau^2_h}{(k+1)^2} \quad \forall k\geq 1,
\end{equation}
where $\Phi(\cdot)$ is the objective function of the dual problem {\rm(\ref{eqn:discretized matrix-vector dual problem})} and
\begin{eqnarray}
  &&\tau^2_h=\frac{1}{2}\|{\bf z}^0- {\bf z}^*\|_{\mathcal{S}^2_h}^2,\label{convergence factor11} \\
  && \mathcal{S}^2_h=\frac{1}{\alpha}\left(
                                      \begin{array}{ccc}
                                        W_h-M_h&  0& 0\\
                                        0&0& 0\\
                                        0&   0& \gamma M_hW_h^{-1}M_h  \\
                                      \end{array}
                                    \right).\label{convergence factor12}
\end{eqnarray}
\end{proposition}

Thus, compared $\mathcal{S}^2_h$ with $\mathcal{S}^1_h$, it is obvious that $\mathcal{S}^2_h\prec\mathcal{S}^1_h$, which implies Algorithm \ref{algo2:Full inexact ABCD algorithm for (Dh)} converge much faster than Algorithm \ref{algo1:Full inexact ABCD algorithm for (Dh)}. However, to solve the $({\bm \lambda}, {\bf p})$-subproblem, utilizing the sGS decomposition technique in Algorithm \ref{algo1:Full inexact ABCD algorithm for (Dh)} would be much easier than using the SNCG algorithm in Algorithm \ref{algo2:Full inexact ABCD algorithm for (Dh)}. Taking the virtues of two variants of the mABCD method into account, we give some strategies about how to choose them. In consideration of the mesh independence and $O(1/k^2)$ iteration complexity, Algorithm sGS-mABCD (Algorithm \ref{algo1:Full inexact ABCD algorithm for (Dh)}) is used to be a priority. In fact, for most of the problems, Algorithm sGS-mABCD can achieve a high accuracy efficiently. However, for some difficult problems, Algorithm sGS-mABCD may not work. In this case, Algorithm mABCD-SNCG (Algorithm \ref{algo2:Full inexact ABCD algorithm for (Dh)}) can perform much better since it makes use of second-order information and it has less blocks. Thus, we can start with Algorithm sGS-mABCD, and then switch to Algorithm mABCD-SNCG when the convergence speed of the Algorithm sGS-mABCD is deemed to be unsatisfactory.

%On the other hand, if we carefully check Algorithm \ref{algo2:Full inexact ABCD algorithm for (Dh)}, we should point out that the information of the accelerated points $\{\tilde {{\bf p}}^k\}_{k>1}$ is not used in the whole iterative process of Algorithm sGS-mABCD-2. Although in theory such the iterative scheme will not affect the convergence result in Proposition \ref{sGS-imABCD convergence}, the lack of such acceleration information may affect the actual convergence rate of the algorithm in the numerical implementation.

\subsection{\textbf{Primal problem of (\ref{eqn:discretized matrix-vector dual problem})}}\label{sec:4.4}
Although we have shown the convergence behavior and the iteration complexity of the dual problem {\rm(\ref{eqn:discretized matrix-vector dual problem})}, our ultimate goal is looking for optimal control solution and optimal state solution. Thus, this can become a driving force for analysing the primal problem of (\ref{eqn:discretized matrix-vector dual problem}). About the primal problem of (\ref{eqn:discretized matrix-vector dual problem}), we have the following result.
\begin{theorem}\label{primal problem of Dh}
Problem {\rm(\ref{eqn:discretized matrix-vector dual problem})} could be regarded as the dual problem of {\rm(\ref{equ:new approx discretized matrix-vector form})}.
In other words, problem {\rm(\ref{equ:new approx discretized matrix-vector form})} is the primal problem of {\rm(\ref{eqn:discretized matrix-vector dual problem})}.
\begin{equation}\label{equ:new approx discretized matrix-vector form}
\left\{\begin{aligned}
        &\min\limits_{{\bf y,u}}^{}~~ {J}_h({\bf y,u})= \frac{1}{2}\|{\bf y-y_{d}}\|_{M_h}^{2}+\frac{\alpha}{2}\|{\bf u}\|_{M_h}^{2}+ \beta\|M_h {\bf u}\|_{1}\\
        &\ {\rm{s.t.}}\quad K_h{\bf y}=M_h ({\bf u +y_r}),\\
        &\ \quad\quad \ {\bf u}\in[a,b]^{N_h}.%\{v\in\mathbb{R}^n: a\leq v_i\leq b, i=1,...,n\}.
                          \end{aligned} \right.%\tag{$\widehat{\mathrm{P}}_{h}$}
\end{equation}
\proof
{Firstly, by introducing two artificial variables, we can rewrite {\rm{(\ref{equ:new approx discretized matrix-vector form})}} as:
\begin{equation}\label{equ:decoupled new approx discretized matrix-vector form}
\left\{\begin{aligned}
        &\min\limits_{\bf y,u,v,w}^{}~~\bar{J}_h({\bf y,u,v,w})= \frac{1}{2}\|{\bf y-y_{d}}\|_{M_h}^{2}+\frac{\alpha}{2}\|{\bf u}\|_{M_h}^{2}+ \beta\|M_h {\bf v}\|_{1}+\delta_{[a,b]}({\bf w})\\
        &\ {\rm{s.t.}}\qquad\quad K_h{\bf y}=M_h ({\bf u +y_r}),\\
        &\ \qquad \qquad M_h({\bf u-v})=0,\\
        &\ \qquad \qquad M_h({\bf u-w})=0.
                          \end{aligned} \right.%\tag{$\widetilde{\mathrm{P}}_{h}$}
\end{equation}
Considering the Lagrangian function associated with {\rm(\ref{equ:decoupled new approx discretized matrix-vector form})}, we have
\begin{eqnarray}\label{Lagrangian function for dual}
\begin{aligned}
  L({\bf y,u,v,w;p,}\bm \lambda,\bm \mu)&=\frac{1}{2}\|{\bf {y-y_{d}}}\|_{M_h}^{2}+\frac{\alpha}{2}\|{\bf u}\|_{M_h}^{2}+ \beta\|M_h {\bf v}\|_{1}+\delta_{[a,b]}({\bf w})\\
  &\quad+\langle {\bf p},K_h{\bf y}-M_h {\bf(u +y_r)}\rangle+\langle\bm \lambda,M_h({\bf u-v})\rangle+\langle\bm \mu,M_h({\bf u-w})\rangle.
\end{aligned}
\end{eqnarray}
Now, we can derive
\begin{eqnarray*}
\left\{\begin{aligned}
  \inf_{\bf y} L({\bf y,u,v,w;p,}\bm \lambda,\bm \mu) &= \inf_{\bf y}{\frac{1}{2}\|{\bf {y-y_{d}}}\|_{M_h}^{2}}+\langle {\bf p},K_h{\bf y}\rangle =-\frac{1}{2}\|K_h {\bf p}-{M_h}{\bf y_{d}}\|_{M_h^{-1}}^2+\frac{1}{2}\|{\bf y_d}\|^2_{M_h},\\
  \inf_{\bf u} L({\bf y,u,v,w;p,}\bm \lambda,\bm \mu)&=\inf_{\bf u} \frac{\alpha}{2}\|{\bf u}\|_{M_h}^{2}-\langle {\bf p},M_h {\bf u}\rangle+\langle\bm \lambda,M_h{\bf u}\rangle+\langle\bm \mu,M_h{\bf u}\rangle=-\frac{1}{2\alpha}\|\bm \lambda+ \bm \mu- {\bf p}\|_{M_h}^2,\\
   \inf_{\bf v}L({\bf y,u,v,w;p,}\bm \lambda,\bm \mu)&=\inf_{\bf v}-\langle\bm \lambda,M_h{\bf v}\rangle+\beta\|M_h{\bf v}\|_1=-\delta_{[-\beta,\beta]}(\bm \lambda), \\
  \inf_{\bf w}L({\bf y,u,v,w;p,}\bm \lambda,\bm \mu)&=\inf_{\bf w}-\langle\bm \mu,M_h{\bf w}\rangle+\delta_{[a,b]}({\bf w})=- \delta^*_{[a,b]}({M_h}\bm \mu).
\end{aligned}\right.
\end{eqnarray*}
Thus,
\begin{eqnarray*}
  \min_{{\bf y,u,v,w}}L({\bf y,u,v,w;p},\bm \lambda,\bm \mu)&=-\frac{1}{2}\|K_h {\bf p}-{M_h}{\bf y_{d}}\|_{M_h^{-1}}^2- \frac{1}{2\alpha}\|\bm \lambda+ \bm \mu- {\bf p}\|_{M_h}^2-\langle M_h{\bf y_r}, {\bf p}\rangle\\
&-\delta_{[-\beta,\beta]}(\bm \lambda)-\delta^*_{[a,b]}({M_h}\bm \mu)+\frac{1}{2}\|{\bf y_d}\|^2_{M_h},
\end{eqnarray*}
and $\max_{{\bf p},\bm \lambda,\bm \mu}\min_{{\bf y,u,v,w}}L({\bf y,u,v,w;p},\bm \lambda,\bm \mu)$ is an equivalent maximization form of the dual problem {\rm(\ref{eqn:discretized matrix-vector dual problem})}. Moveover, there is no gap between {\rm(\ref{equ:new approx discretized matrix-vector form})} and {\rm(\ref{eqn:discretized matrix-vector dual problem})} due to the strong convexity of problem {\rm (\ref{equ:new approx discretized matrix-vector form})}. Thus, we complete the proof.}
\endproof
\end{theorem}

Since the stiffness matrix $K_h$ is a symmetric positive definite matrix, problem (\ref{equ:new approx discretized matrix-vector form}) can be rewritten as the following reduced form:

\begin{equation}\label{equ:reduced new approx discretized matrix-vector form}
        \min\limits_{\bf u}^{}~~\widehat{J}_h({\bf u})= \frac{1}{2}\|K_h^{-1}M_h{\bf (u+y_r)-y_{d}}\|_{M_h}^{2}+\frac{\alpha}{2}\|{\bf u}\|_{M_h}^{2}+ \beta\|M_h {\bf u}\|_{1}+\delta_{[a,b]}({\bf u}).\\
%\tag{$\widehat{\mathrm{RP}}_{h}$}
\end{equation}
Thus, about the iteration complexity of the primal problem {\rm(\ref{equ:new approx discretized matrix-vector form})} of (\ref{eqn:discretized matrix-vector dual problem}), we have the following results.

\begin{proposition}{\rm{\textbf{\cite[Theorem 10]{meshindependencemABCD}}}}\label{iteration complexity of the primal problem}
Let $\{{\bf z}^k\}:=\{(\bm \lambda^k,{\bf p}^k,\bm \mu^k)\}$ be the sequence generated by the Algorithm {\rm \ref{algo1:Full inexact ABCD algorithm for (Dh)}}, ${\bf u}^k=({\bf p}^k-\bm \lambda^k-\bm \mu^k)/\alpha$ and denote  $\hat{{\bf u}}^k={\rm \Pi}_{[a,b]}({\bf u}^k)$, then
\begin{equation}\label{difference of primal and dual obejective function}
\widehat{J}_h(\hat{\bf u}^k)-\widehat{J}_h({\bf u}^*)\leq\frac{C_1}{1+k} \quad \forall k\geq 1,
\end{equation}
where $\widehat{J}_h$ is the objective function of problem {\rm(\ref{equ:reduced new approx discretized matrix-vector form})} and ${\bf u}^*$ is the unique optimal solution of problem {\rm(\ref{equ:reduced new approx discretized matrix-vector form})}, moreover, we have
\begin{equation}\label{the convergence of the primal sequence}
\|\hat{\bf u}^k-{\bf u}^*\|\leq \frac{C_2}{\sqrt{k+1}} \quad \forall k\geq 1,
\end{equation}
where $C_1$ and $C_2$ are two constants.
\end{proposition}

\section{\textbf{Error estimates}}\label{sec:5}
Based on the Theorem \ref{primal problem of Dh}, for any $u_h=\sum\limits^{N_h}_{i=1}u_i\phi_i$ and ${\bf u}=(u_1,u_2,...,u_{N_h})$, now let us define the new approximation of the $L^1$-norm by
\begin{equation}\label{new approx discretization21}
\|u_h\|_{\widetilde{L}^1_h(\Omega)}=\sum_{j=1}^{N_h}|\int_{\Omega_h}\sum\limits_{i=1}^{N_h}u_i\phi_i(x)\phi_j(x)\mathrm{d}x|=\|M_h\mathbf{u}\|_1,
\end{equation}
which can be regarded as a generalized weighted $l^1$-norm of the coefficients of $u_h$, thus it is a norm on $U_h$. Moreover, we can rewrite {\rm(\ref{equ:new approx discretized matrix-vector form})} as the following discretized function form:
\begin{equation}\label{eqn:new approx discretized problems}
  \left\{ \begin{aligned}
        &\min \limits_{(y_h,u_h)\in Y_h\times U_h}^{}J_h(y_h,u_h)=\frac{1}{2}\|y_h-y_d\|_{L^2(\Omega_h)}^{2}+\frac{\alpha}{2}\|u_h\|_{L^2(\Omega_h)}^{2}+ \beta\|u_h\|_{\widetilde{L}_h^1(\Omega_h)}\\
        &\qquad\quad {\rm{s.t.}}\qquad a(y_h, v_h)=\int_{\Omega}(u_h+y_r)v_h{\rm{d}}x, \qquad  \forall v_h\in Y_h,  \\
          &\qquad \qquad \qquad~  u_h\in U_{ad,h},
                          \end{aligned} \right.\tag{$\widehat{\mathrm{P}}_{h}$}
 \end{equation}
%where $\|\cdot\|_{\widetilde{L}_h^1(\Omega_h)}$ is defined as $\|u_h\|_{\widetilde{L}_h^1(\Omega_h)}=\|M_hu\|_1$ and $u_h=\sum_{i=1}^{N_h}{u_i\phi_i(x)}$.

In this section, we will accomplish the error estimates for the discretized problem (\ref{eqn:new approx discretized problems}).

\subsection{\textbf{Analysis of the approximate $L^1$-norm}}\label{sec:5.1}
Since $\|M_hu\|_1$ can be regarded as an approximation of $\|u_h\|_{L^1(\Omega)}=\int_{\Omega_h}|\sum_{i=1}^{N_h}{u_i\phi_i(x)}|~\mathrm{d}x$, it is necessarily required to analyse the finite element error. First, we will analyse the relationship between $\|M_hu\|_1$ and $\|u_h\|_{L^1(\Omega)}$ . For the analyse further below, let us first introduce a quasi-interpolation operator $\Pi_h:L^1(\Omega_h)\rightarrow U_h$ which provides interpolation estimates. For an arbitrary $w\in L^1(\Omega)$, the operator $\Pi_h$ is constructed as follows:
\begin{equation}\label{equ:quasi-interpolation}
   \Pi_hw=\sum\limits_{i=1}^{N_h}\pi_i(w)\phi_i(x), \quad \pi_i(w)=\frac{\int_{\Omega_h}w(x)\phi_i(x){\rm{d}}x}{\int_{\Omega_h}\phi_i(x){\rm{d}}x}.
\end{equation}
Moreover, since the upper and lower bounds $a$ and $b$ are constants, we have that
\begin{equation}\label{quasi-interpolation property}
 w\in U_{ad} \Rightarrow \Pi_hw \in U_{ad,h} \quad {\rm for\  all}\ w\in L^1(\Omega).
\end{equation}
Based on the assumption on the mesh and the control discretization, we extend $\Pi_hw$ to $\Omega$ by taking $\Pi_hw=w$ for every $x\in\Omega\backslash {\Omega}_h$, and we have the following estimates of the interpolation error. For the detailed proofs, we refer to \cite{Cars}.
\begin{lemma}\label{eqn:lemma3}
There is a constant $c_{\pi}$ independent of $h$ such that
\begin{equation}\label{error-estimates-quasi-interpolation}
  h\|z-\Pi_hz\|_{L^2(\Omega)}+\|z-\Pi_hz\|_{H^{-1}(\Omega)}\leq c_{\pi}h^2\|z\|_{H^1(\Omega)}
\end{equation}
holds for all $z\in H^1(\Omega)$.
\end{lemma}

Then, from the definition of $\Pi_h$ and Lemma \ref{eqn:lemma3}, we have the following results.
\begin{proposition}\label{eqn:another martix properties}
$\forall$ ${\bf z}=(z_1,z_2,...,z_{N_h})\in \mathbb{R}^{N_h}$ and $z_h(x)=\sum\limits_{i=1}^{N_h}z_i\phi_i$, there exists a constant $C$, such that the following inequalities hold
\begin{eqnarray}
  \label{Minequality1}&&\|\Pi_h z_h\|_{L^1(\Omega)}\leq\|z_h\|_{\widetilde{L}^1_h(\Omega)}\leq\|z_h\|_{L^1(\Omega)},\\%\int_{\Omega_h}|\sum_{i=1}^{N_h}{z_i\phi_i(x)}|~\mathrm{d}x,\\
  \label{Minequality2}&&\|z_h\|_{L^1(\Omega)}-\|z_h\|_{\widetilde{L}^1_h(\Omega)}\leq C\|z_h\|_{H^1(\Omega)}h^2,
  %\int_{\Omega_h}|\sum_{i=1}^n{z_i\phi_i(x)}|~\mathrm{d}x-\|M_hz\|_1=C\|z_h\|_{H^1(\Omega)}h^2,
\end{eqnarray}
where $C$ is independent of $h$.
\begin{proof}
Since
\begin{equation*}
\begin{aligned}
 \|z_h\|_{\widetilde{L}^1_h(\Omega)}=\|M_h{\bf z}\|_1&=\sum_{i=1}^{N_h}|\sum_{j=1}^{N_h}\int_{\Omega_h}{\phi_i(x)\phi_j(x)z_j}~\mathrm{d}x|\\
 &\leq\sum_{i=1}^{N_h}\int_{\Omega_h}|\sum_{j=1}^{N_h}{\phi_j(x)z_j}|\phi_i(x)~\mathrm{d}x\\
 &=\int_{\Omega_h}|\sum_{i=1}^{N_h}{z_i\phi_i(x)}|~\mathrm{d}x=\|z_h\|_{L^1(\Omega)},
\end{aligned}
\end{equation*}
where the last equality is due to $\sum\limits_{i=1}^{N_h}\phi_i=1$. Furthermore,
\begin{equation*}
\begin{aligned}
  \|\Pi_h z_h\|_{L^1(\Omega)}&=\int_{\Omega_h}|\sum_{j=1}^{N_h}\frac{\int_{\Omega_h}\sum\limits_{i=1}^{N_h}z_i\phi_i(x)\phi_j(x)\mathrm{d}x}{\int_{\Omega_h}\phi_j(x)\mathrm{d}x}\phi_j(x)|\mathrm{d}x\\
  &\leq\sum_{j=1}^{N_h}\frac{|\int_{\Omega_h}\sum\limits_{i=1}^{N_h}z_i\phi_i(x)\phi_j(x)\mathrm{d}x|}{\int_{\Omega_h}\phi_j(x)\mathrm{d}x}\int_{\Omega_h}\phi_j(x)\mathrm{d}x\\
  &=\|M_h{\bf z}\|_1=\|z_h\|_{\widetilde{L}^1_h(\Omega)}.
\end{aligned}
\end{equation*}
At last, according to $z_h\in H^1(\Omega_h)$ and Lemma {\rm\ref{eqn:lemma3}}, we have
\begin{eqnarray*}
  \|z_h\|_{L^1(\Omega)}-\|z_h\|_{\widetilde{L}^1_h(\Omega)}
  &\leq&\|z_h\|_{L^1(\Omega)}-\|\Pi_h z_h\|_{L^1(\Omega)}\\
  &\leq&\|\Pi_h z_h-z_h\|_{L^1(\Omega)}\\
  &\leq& c_{\Omega}\|\Pi_h z_h- z_h\|_{H^{-1}(\Omega)}\\
  &\leq&c_{\Omega}c_{\pi}\|z_h\|_{H^1(\Omega)}h^2\\
  &=&C\|z_h\|_{H^1(\Omega)}h^2.
\end{eqnarray*}
\end{proof}
\end{proposition}

Thus, from Proposition \ref{eqn:another martix properties}, it is reasonable to consider {\rm(\ref{eqn:new approx discretized problems})} as a discretization of problem {\rm(\ref{eqn:orginal problems})}. Finally, compared (\ref{Minequality2}) with (\ref{Winequality4}), it is obvious that utilizing $\|u_h\|_{\widetilde{L}_h^1(\Omega_h)}$ to approximate $\|u_h\|_{L^1(\Omega_h)}$ is better than using $\|u_h\|_{L_h^1(\Omega_h)}$. Hence, from the point of view of the approximation order, equivalently solving the dual problem {\rm(\ref{eqn:discretized matrix-vector dual problem})} of (\ref{eqn:new approx discretized problems}) is superior to solving (\ref{eqn:approx discretized problems}). %Actually, in \cite{mABCDSOPT}, from their numerical results, the authors have already shown that solving the dual problem (\ref{eqn:discretized matrix-vector dual problem}) could get better error results than that from solving (\ref{eqn:approx discretized problems}). %Thus, in consideration of this phenomenon, another key issue should be considered is how measures of the solution accuracy of (\ref{eqn:new approx discretized problems}) vary with the level of discretized approximation. Such questions come under the category of the finite element error estimates.

\subsection{\textbf{Finite element error estimates}}\label{sec:5.2}
Now, let us return to the a-priori error analysis. Analogous to the continuous problem (\ref{eqn:orginal problems}), the discretized problem (\ref{eqn:new approx discretized problems}) is also a strictly convex problem, which is uniquely solvable. We derive the following first-order optimality conditions, which are necessary and sufficient for the optimal solution of (\ref{eqn:new approx discretized problems}).
\begin{theorem}[{\rm \textbf{Discretized first-order optimality condition}}]
$(y^*_h, u^*_h)$ is the optimal solution of {\rm(\ref{eqn:new approx discretized problems})} if and only if there exists an adjoint state $p^*_h$, such that the following conditions are satisfied
\begin{subequations}\label{eqn:DKKT}
\begin{eqnarray}
&&y^*_h=\mathcal{S}_h (u^*_h+y_r)\label{eqn1:DKKT}, \\
&&p^*_h=\mathcal{S}_h^*(y_d-y^*_h)\label{eqn2:DKKT},\\
&&{\left\langle\alpha u^*_h-p^*_h,u_h-u^*_h\right\rangle_{L^{2}(\Omega_h)}+\beta\left(\|u_h\|_{\widetilde{L}^{1}_h(\Omega_h)}-\|u^*_h\|_{\widetilde{L}^{1}_h(\Omega_h)}\right)}\geq0 \label{eqn3:DKKT}\quad\forall u_h \in U_{ad,h}.
\end{eqnarray}
\end{subequations}
\end{theorem}

Before we analyse the finite element error estimates, let us introduce an important inequality, which is
called scaling argument. For two Banach spaces $B_0$, $B_1$, the continuous embedding $B_1\hookrightarrow B_0$ implies that there exists a constant $C>0$ such that
\begin{equation*}
  \|v\|_{B_0}\leq C\|v\|_{B_1}, \quad \forall v\in B_1.
\end{equation*}
The inequality in the reserve way $\|v\|_{B_1}\leq C\|v\|_{B_0}$ may not true. However, considering finite
element spaces $V_h \subset B_i, i=0,1$ endowed with two norms, since the dimension
of $V_h$ is finite and all the norms of finite dimensional spaces are equivalent, the above inverse inequality
will be true for all $v\in V_h$. This result is shown as below.

\begin{lemma}\label{lemma typical inverse inequalities}
For all $v_h\in U_h$, $0<h<h_0$, there exists a constant $C$ such that
\begin{equation}\label{typical inverse inequalities}
  \|v_h\|_{L^2(\Omega)}\leq\|v_h\|_{H^1(\Omega)}\leq Ch^{-1}\|v_h\|_{L^2(\Omega)}.
\end{equation}
\begin{proof}
Due to the continuous embedding $H^1(\Omega)\hookrightarrow L^2(\Omega)$, it is obvious that we have $\|v_h\|_{L^2(\Omega)}\leq\|v_h\|_{H^1(\Omega)}$. In addition, by $v_h=\sum_{i=1}^{N_h}v_i\phi(x)$, ${\bf v}=(v_1,v_2,...,v_{N_h})$ and Lemma {\rm\ref{spectral property}}, we have
\begin{equation*}
  \begin{aligned}
  \frac{\|v_h\|_{H^1(\Omega)}}{\|v_h\|_{L^2(\Omega)}}&=\sqrt{\frac{{\bf v}^{T}M_h\mathbf{v}+\mathbf{v}^{T}K_h\mathbf{v}}{\mathbf{v}^{T}M_h\mathbf{v}}}
                                                     \leq\left\{\begin{aligned}
                                                            \sqrt{1+\frac{d_2}{c_1h^2}}\quad {\rm for}~n=2  \\
                                                            \sqrt{1+\frac{d_2h}{c_1h^3}}\quad {\rm for}~n=3  \\
                                                          \end{aligned}\right.
                                                     \leq\sqrt{\frac{c_1h_0^2+d_2}{c_1h^2}}
                                                     =Ch^{-1}.
  \end{aligned}
\end{equation*}
Thus, we complete the proof.
\end{proof}
\end{lemma}

Now, let us start to derive error estimation in terms of the mesh size with the help of the
variational inequalities (\ref{eqn3:KKT}) and (\ref{eqn3:DKKT}).
\begin{theorem}\label{theorem:error2}
Let $(y^*, u^*)$ be the optimal solution of problem {\rm(\ref{eqn:orginal problems})}, and $(y^*_h, u^*_h)$ be the optimal solution of problem {\rm(\ref{eqn:new approx discretized problems})}. For any $h>0$ small enough and $\alpha_0,\beta_0>0$, there is a constant $C$ such that for all $0<\alpha\leq\alpha_0$ and $0<\beta\leq\beta_0$,
\begin{eqnarray*}
  &&\frac{\alpha}{4}\|u^*-u^*_h\|^2_{L^2(\Omega)}+\frac{1}{2}\|y^*-y^*_h\|^2_{L^2(\Omega)}\\
  &&\leq C(\alpha h^2+h^2+\alpha^{-1}\beta^2 h^2+\alpha^{-1}\beta h^2+\alpha^{-1} h^2+\alpha^{-1}h^3+h^4+\alpha^{-1}h^4+\alpha^{-2}h^4),
\end{eqnarray*} where $C$ is a constant independent of $h$, $\alpha$ and $\beta$.
\begin{proof}
Due to the optimality of $u^*$ and $u^*_h$, $u^*$ and $u^*_h$ satisfy {\rm(\ref{eqn3:KKT})} and {\rm(\ref{eqn3:DKKT})}, respectively. From $U_{ad,h}\subset U_{ad}$, thus the function $u^*_h$
is feasible for the continuous problem, i.e. $u^*_h\in U_{ad}$, and can be used as test
function in the variational inequality {\rm(\ref{eqn3:KKT})}. On the other hand, it would be nice if we could use $u^*$ as a test function in the variational
inequality {\rm(\ref{eqn3:DKKT})}, which characterizes $u^*_h$. However, in general the function $u^*$ does
not belong to $U_{ad,h}$ and cannot be utilized as test function. To overcome this difficulty, let us introduce an approximation $\widetilde{u}^*_h:=\Pi_hu^*\in U_{ad,h}$, which is suitable as test function in {\rm(\ref{eqn3:DKKT})}.

Now, let us use the test function $u^*_h$ in {\rm(\ref{eqn3:KKT})} and the test function $\widetilde{u}^*_h$ in {\rm(\ref{eqn3:DKKT})}, thus we have
\begin{eqnarray}
\label{eqn:error1}&&{\left\langle\alpha u^*-p^*,u^*_h-u^*\right\rangle_{L^2(\Omega)}+\beta\left(\|u^*_h\|_{L^1(\Omega)}-\|u^*\|_{L^1(\Omega)}\right)}\geq0, \\
\label{eqn:error2}&&{\left\langle\alpha u^*_h-p^*_h,\widetilde{u}^*_h-u^*_h\right\rangle_{L^{2}(\Omega_h)}+\beta\left(\|\widetilde{u}^*_h\|_{\widetilde{L}^{1}_h(\Omega_h)}-\|u^*_h\|_{\widetilde{L}^{1}_h(\Omega_h)}\right)}\geq0.
\end{eqnarray}
Because $u^*_h=0$ on $\bar\Omega\backslash {\Omega}_h$, the integrals over $\Omega$ can be replaced by integrals over $\Omega_h$ in {\rm(\ref{eqn:error1})}, and it can be rewritten as
\begin{equation}\label{eqn:error3}
  \begin{aligned}
  &{\left\langle\alpha u^*-p^*,u^*-u^*_h\right\rangle_{L^2(\Omega_h)}+\beta\left(\|u^*\|_{L^1(\Omega_h)}-\|u^*_h\|_{L^1(\Omega_h)}\right)}\\
  &\leq\left\langle p^*-\alpha u^*,u^*\right\rangle_{L^2(\Omega\backslash {\Omega}_h)}-\beta\|u^*\|_{L^1(\Omega\backslash {\Omega}_h)} \\
  &\leq \langle p^*,u^*\rangle_{L^2(\Omega\backslash {\Omega}_h)}\leq C_1h^2,
  \end{aligned}
\end{equation}
where the last inequality follows from the boundedness of $p^*$ and $u^*$ and the assumption $|\Omega\backslash {\Omega}_h|\leq c h^2$.
In addition, {\rm(\ref{eqn:error2})} can be rewritten as
\begin{eqnarray}
\label{eqn:error4}&&{\left\langle p^*_h- \alpha u^*_h,u^*-u^*_h+\widetilde{u}^*_h-u^*\right\rangle_{L^{2}(\Omega_h)}+\beta\left(\|u^*_h\|_{\widetilde{L}^{1}_h(\Omega_h)}-\|\widetilde{u}^*_h\|_{\widetilde{L}^{1}_h(\Omega_h)}\right)}\leq0.
\end{eqnarray}
Adding up and rearranging {\rm(\ref{eqn:error3})} and {\rm(\ref{eqn:error4})}, we obtain
\begin{equation}\label{eqn:error estimat1}
  \begin{aligned}
  \alpha\|u^*-u^*_h\|^2_{L^2(\Omega_h)}\leq& \langle p^*-p^*_h,u^*-u^*_h\rangle_{L^2(\Omega_h)}+\langle \alpha u^*_h -p^*_h,\widetilde{u}^*_h-u^*\rangle_{L^2(\Omega_h)}\\
&+\beta\left(\|u^*_h\|_{L^1(\Omega_h)}-\|u^*_h\|_{\widetilde{L}_h^1(\Omega_h)}+\|\widetilde{u}^*_h\|_{\widetilde{L}^{1}_h(\Omega_h)}-\|u^*\|_{L^{1}(\Omega_h)}\right)+C_1h^2\\
\leq&
\underbrace{\left\langle\alpha u^*-p^*,\widetilde{u}^*_h-u^*\right\rangle_{L^2(\Omega_h)}}
_{I_1}
+
\underbrace{\alpha\left\langle u^*_h-u^*,\widetilde{u}^*_h-u^*\right\rangle_{L^2(\Omega_h)}}_{I_2}
+\underbrace{\left\langle p^*-p^*_h,\widetilde{u}^*_h-u_h^*\right\rangle_{L^2(\Omega_h)}}_{I_3}\\
&+\underbrace{\beta\left(\|u^*_h\|_{L^1(\Omega_h)}-\|u^*_h\|_{\widetilde{L}_h^1(\Omega_h)}+\|\widetilde{u}^*_h\|_{\widetilde{L}^{1}_h(\Omega_h)}-\|u^*\|_{L^{1}(\Omega_h)}\right)}
_{I_4}
+C_1h^2,
\end{aligned}
\end{equation}
First, let us estimate the fourth term $I_4$. From the definition of $\widetilde{u}^*_h=\Pi_h(u^*)$ and the non-negativity and partition of unity of the nodal basis functions, we get
\begin{equation}\label{equ:equation about l^1 norm}
\begin{aligned}
  \|\widetilde{u}^*_h\|_{\widetilde{L}^{1}_h(\Omega_h)}=\|\Pi_h(u^*)\|_{\widetilde{L}^{1}_h(\Omega_h)}
&=\sum\limits_{i=1}^{N_h}\left|\sum\limits_{j=1}^{N_h}{\int_{\Omega_h}\phi_i(x)\phi_j(x)\frac{\int_{\Omega_h}u^*\phi_j(x){\rm{d}}x}{\int_{\Omega_h}\phi_j(x){\rm{d}}x}{\rm{d}}x}\right|\\
&\leq\sum\limits_{i=1}^{N_h}\int_{\Omega_h}\sum\limits_{j=1}^{N_h}\left|\frac{\int_{\Omega_h}u^*\phi_j(x){\rm{d}}x}{\int_{\Omega_h}\phi_j(x){\rm{d}}x}\right|{\phi_i(x)\phi_j(x){\rm{d}}x}
=\sum\limits_{j=1}^{N_h}\left|\int_{\Omega_h}u^*\phi_j(x){\rm{d}}x\right|\\
&\leq\sum\limits_{j=1}^{N_h}\int_{\Omega_h}\left|u^*\right|\phi_j(x){\rm{d}}x
=\|u^*\|_{L^1(\Omega_h)}.
\end{aligned}
\end{equation}
Thus, it suffices to employ Proposition {\rm\ref{eqn:another martix properties}} and Lemma {\rm\ref{interpolation error estimate}}, in conjunction with Lemma {\rm\ref{lemma typical inverse inequalities}},
we obtain
\begin{equation}\label{I_4 term}
\begin{aligned}
\|u^*_h\|_{L^1(\Omega_h)}-\|u^*_h\|_{\widetilde{L}_h^1(\Omega_h)}&\leq C_2\|u^*_h\|_{H^1(\Omega)}h^2\\
&\leq C_2 h^2(\|u^*_h-u^*\|_{H^1(\Omega)}+\|u^*\|_{H^1(\Omega)})\\
&\leq C_2 h^2(\|u^*_h-I_hu^*\|_{H^1(\Omega)}+\|u^*-I_hu^*\|_{H^1(\Omega)}+\|u^*\|_{H^1(\Omega)})\\
&\leq C_2C_3 h(\|u^*_h-I_hu^*\|_{L^2(\Omega)})+C_2(c_I+1)h^2\|u^*\|_{H^1(\Omega)}\\
&\leq C_4 h(\|u^*_h-u^*\|_{L^2(\Omega)}+\|u^*-I_hu^*\|_{L^2(\Omega)})+C_5h^2\|u^*\|_{H^1(\Omega)}\\
&\leq C_4 h\|u^*_h-u^*\|_{L^2(\Omega)}+C_4c_Ih^2\|u^*\|_{H^1(\Omega)}+C_5h^2\|u^*\|_{H^1(\Omega)}\\
&\leq C_6\beta\alpha^{-1}h^2+ \frac{\alpha}{8\beta}\|u^*-u^*_h\|^2_{L^2(\Omega)}+C_7h^2\|u^*\|_{H^1(\Omega)}\\
%&\leq C_3 h^2(C'\inf_{v_h\in U_{h}}\|u^*-v_h\|_{H^1(\Omega)}+\|u^*\|_{H^1(\Omega)})
%&\leq C h^2(C'\|u^*-I_hu^*\|_{H^1(\Omega)}+\|u^*\|_{H^1(\Omega)})\\
%&\leq C h^2(C'c_I\|u^*\|_{H^1(\Omega)}+\|u^*\|_{H^1(\Omega)})\\
%&=C_2\|u^*\|_{H^1(\Omega)}h^2
\end{aligned}
\end{equation}
From the regularity of the optimal control $u^*$ and optimal adjoint state $p^*$, i.e. $u^*\in H^1(\Omega)$ and $p^*\in H^1(\Omega)$, and {\rm(\ref{eqn3.1:KKT})}, we know that
\begin{equation}\label{eqn:exact function estimats1}
  \|u\|_{H^1(\Omega)}\leq \frac{1}{\alpha}\|p\|_{H^1(\Omega)}+\left(\frac{\beta}{\alpha}+|a|+b\right)\mathcal{M}(\Omega)^{\frac{1}{2}},
\end{equation}
where $\mathcal{M}(\Omega)$ denotes the measure of $\Omega$. Obviously, the $H^1$-norm of $u^*$ depends on $\alpha$. On the other hand, due to the control constraints, the $H^1$-norm of $p^*$ is bounded independently of $\alpha$. Thus we have
\begin{equation*}
  \|\alpha u^*-p^*\|_{H^1(\Omega)}\leq2\|p\|_{H^1(\Omega)}+(\beta+\alpha |a|+\alpha b)\mathcal{M}(\Omega)^{\frac{1}{2}}.
\end{equation*}
Moreover, due to the boundedness of the state $y^*$, the adjoint state $p^*$ and the operator $\mathcal{S}$, we can choose a large enough constant $L>0$ independent of $\alpha$, $\beta$ and $h$ and two constants $\alpha_0$ and $\beta_0$ such that for all $0<\alpha\leq\alpha_0$, $0<\beta\leq\beta_0$ and $h>0$, the following inequation holds:
\begin{equation}\label{eqn:exact function estimats2}
  2\|p^*\|_{H^1(\Omega)}+(\beta+\alpha |a|+\alpha b+|a|+b)\mathcal{M}(\Omega)^{\frac{1}{2}}+\|y^*-y_d\|_{L^2(\Omega)}+\|y_r\|_{L^2(\Omega)}+\|\mathcal{S}\|_{\mathcal{L}(H^{-1},L^2)}\leq L.
\end{equation}
Thus, we have
\begin{equation}\label{estimatesI4}
  I_4\leq C_6\beta^2\alpha^{-1}h^2+C_7L\beta\alpha^{-1} h^2+\frac{\alpha}{8}\|u^*-u^*_h\|^2_{L^2(\Omega)}.
\end{equation}

Next, in order to further estimate {\rm(\ref{eqn:error estimat1})}, we will discuss each of these items from $I_1$ to $I_3$ in turn. For the terms $I_1$, from the regularity of the optimal control $u^*$ and optimal adjoint state $p^*$, i.e. $u^*\in H^1(\Omega)$ and $p^*\in H^1(\Omega)$, and Lemma {\rm\ref{eqn:lemma3}}, we get
\begin{equation}\label{estimatesI1}
\begin{aligned}
I_1=&\left\langle\alpha u^*-p^*,\widetilde{u}^*_h-u^*\right\rangle_{L^2(\Omega_h)}\\
\leq& \|\alpha u^*-p^*\|_{H^1(\Omega_h)}\|\widetilde{u}^*_h-u^*\|_{H^{-1}(\Omega_h)}\\
\leq&Lc_{\pi}\|u\|_{H^1(\Omega)}h^2\\
\leq&\alpha^{-1}c_{\pi}L^2h^2
\end{aligned}
\end{equation}
%where the last inequality is due to (\ref{eqn:exact function estimats2}).
For the term $I_2$, we have
\begin{equation}\label{estimatesI2}
\begin{aligned}
I_2=&\alpha\left\langle u^*_h-u^*,\widetilde{u}^*_h-u^*\right\rangle_{L^2(\Omega_h)}\\
\leq& \frac{\alpha}{2}\|u^*_h-u^*\|^2_{L^2(\Omega_h)}+\frac{\alpha}{2}\|\widetilde{u}^*_h-u^*\|^2_{L^{2}(\Omega_h)}\\
\leq&\frac{\alpha}{2}\|u^*_h-u^*\|^2_{L^2(\Omega_h)}+\frac{\alpha}{2}c^2_{\pi}\|u\|_{H^1(\Omega)}^2 h^2\\
\leq&\frac{\alpha}{2}\|u^*_h-u^*\|^2_{L^2(\Omega_h)}+\alpha^{-1}L^2 c^2_{\pi} h^2
\end{aligned}
\end{equation}
For the term $I_3$, let $\tilde{p}^*_h=\mathcal{S}^*_h(y_d-y^*)$, we have
\begin{equation}\label{estimatesI3}
\begin{aligned}
  I_3=&\left\langle p^*-p^*_h,\widetilde{u}^*_h-u_h^*\right\rangle_{L^2(\Omega_h)}\\
=&\langle \tilde{p}^*_h-p^*_h+p^*-\tilde{p}^*_h,\widetilde{u}^*_h-u_h^*\rangle_{L^2(\Omega_h)} \\
=&\left\langle\mathcal{S}^*_h(y^*_h-y^*),\widetilde{u}^*_h-u_h^*\right\rangle_{L^2(\Omega_h)}+\left\langle(\mathcal{S}^*_h-\mathcal{S}^*)(y^*-y_d),\widetilde{u}^*_h-u_h^*\right\rangle_{L^2(\Omega_h)}\\
=&\left\langle y^*_h-y^*,\mathcal{S}_h(\widetilde{u}^*_h-u_h^*)\right\rangle_{L^2(\Omega_h)}+\left\langle y^*-y_d,(\mathcal{S}_h-\mathcal{S})(\widetilde{u}^*_h-u_h^*)\right\rangle_{L^2(\Omega_h)}\\
=&\langle y^*_h-y^*,\mathcal{S}_h\widetilde{u}^*_h-\mathcal{S}u^*+\mathcal{S}u^*-\mathcal{S}_hu^*_h\rangle_{L^2(\Omega_h)}+\langle y^*-y_d,(\mathcal{S}_h-\mathcal{S})(\widetilde{u}^*_h-u^*_h)\rangle_{L^2(\Omega_h)}\\
     =&-\|y^*-y^*_h\|^2_{L^2(\Omega_h)}+
 \underbrace{\langle y^*_h-y^*,(\mathcal{S}_h-\mathcal{S})(\widetilde{u}^*_h+y_r)-\mathcal{S}(u^*-\widetilde{u}^*_h)\rangle_{L^2(\Omega_h)}}
 _{I_5}
\\
 &+\underbrace{\left\langle y^*-y_d,(\mathcal{S}_h-\mathcal{S})(\widetilde{u}^*_h-u_h^*)\right\rangle_{L^2(\Omega_h)}}_{I_6}
\end{aligned}
\end{equation}
For terms $I_5$ and $I_6$, using H$\mathrm{\ddot{o}}$lder's inequality, Lemma {\rm\ref{eqn:lemma1}} and Lemma {\rm\ref{eqn:lemma3}}, we have
\begin{equation}\label{estimatesI5}
  \begin{aligned}
  I_5&\leq\frac{1}{2}\|y^*-y^*_h\|^2_{L^2(\Omega_h)}+2\|\mathcal{S}_h-\mathcal{S}\|^2_{\mathcal{L}(L^{2},L^2)}
  (\|\widetilde{u}^*_h\|^2_{L^2(\Omega_h)}+\|y_r\|^2_{L^2(\Omega_h)})
+\|\mathcal{S}\|_{\mathcal{L}(H^{-1},L^2)}\|u^*-\widetilde{u}^*_h\|^2_{H^{-1}(\Omega_h)}\\
  &\leq\frac{1}{2}\|y^*-y_h\|^2_{L^2(\Omega_h)}+2c_1^2h^4(\sup_{u_h\in U_{ad,h}}\|u_h\|^2_{L^2(\Omega_h)}+\|y_r\|^2_{L^2(\Omega_h)})+c^2_{\pi}L\|u^*\|^2_{H^1(\Omega_h)}h^4\\
  &\leq\frac{1}{2}\|y^*-y_h\|^2_{L^2(\Omega_h)}+2c_1^2L^2h^4+c_{\pi}^2L^3\alpha^{-2}h^4,
\end{aligned}
\end{equation}
and
\begin{equation}\label{estimatesI6}
  \begin{aligned}
  I_6&\leq \|y^*-y_d\|_{L^2(\Omega_h)}\|\mathcal{S}_h-\mathcal{S}\|_{\mathcal{L}(L^{2},L^2)}(\|\widetilde{u}^*_h-u^*\|_{L^2(\Omega_h)}+\|u^*-u^*_h\|_{L^2(\Omega_h)})\\
  &\leq c_1Lh^2(c_{\pi}L\alpha^{-1}h+\|u^*-u^*_h\|_{L^2(\Omega_h)})\\
  &= c_1Lh^2\|u^*-u^*_h\|_{L^2(\Omega_h)}+c_1c_{\pi}\alpha^{-1}L^2h^3\\
  &\leq \frac{\alpha}{4}\|u^*-u^*_h\|^2_{L^2(\Omega_h)}+c_1c_{\pi}\alpha^{-1}L^2h^3+4c_1^2L^2\alpha^{-1}h^4.
\end{aligned}
\end{equation}
Consequently, substituting {\rm(\ref{estimatesI4})}, {\rm(\ref{estimatesI1})}, {\rm(\ref{estimatesI2})}, {\rm(\ref{estimatesI3})}, {\rm(\ref{estimatesI5})} and {\rm(\ref{estimatesI6})} into {\rm(\ref{eqn:error estimat1})} and rearranging, we get
\begin{equation}\label{error estimat2}
\begin{aligned}
  &\frac{\alpha}{8}\|u^*-u^*_h\|^2_{L^2(\Omega_h)}+\frac{1}{2}\|y^*-y^*_h\|^2_{L^2(\Omega_h)}\\ &\leq(C_1+C_6\beta^2\alpha^{-1} +C_7L\beta\alpha^{-1} +\alpha^{-1}c_{\pi}L^2+\alpha^{-1}L^2 c^2_{\pi} )h^2+c_1c_{\pi}\alpha^{-1}L^2h^3\\
  &\quad+(2c_1^2L^2+c_{\pi}^2L^3\alpha^{-2}+4c_1^2L^2\alpha^{-1})h^4\\
  &=C(\alpha h^2+h^2+\alpha^{-1}\beta^2 h^2+\alpha^{-1}\beta h^2+\alpha^{-1} h^2+\alpha^{-1}h^3+h^4+\alpha^{-1}h^4+\alpha^{-2}h^4)
  \end{aligned}
\end{equation}
where $C>0$ is a properly chosen constant.
Thus, the proof is completed.
\end{proof}
\end{theorem}

\begin{corollary}\label{corollary:error1}
Let $(y^*, u^*)$ be the optimal solution of problem {\rm(\ref{eqn:orginal problems})}, and $(y^*_h, u^*_h)$ be the optimal solution of problem {\rm(\ref{eqn:new approx discretized problems})}. For every $h_0>0$, $\alpha_0>0$ and $\beta_0>0$, there is a constant $C>0$ such that for all $0<\alpha\leq\alpha_0$, $0<\beta\leq\beta_0$, $0<h\leq h_0$ it holds
\begin{eqnarray*}
  \|u-u_h\|_{L^2(\Omega)}\leq C(\alpha^{-1}h+\alpha^{-\frac{3}{2}}h^2),
\end{eqnarray*}
where $C$ is a constant independent of $h$, $\alpha$ and $\beta$.
\end{corollary}

\section{Concluding remarks}\label{sec:6}
In this paper, instead of directly solving the primal optimal control problem with $L^1$ control cost, we introduce a duality-based approach and an accelerated block coordinate descent (ABCD) method to solve this type of problem via its dual. Some convergence results for the dual problem are presented. In consideration of our ultimate goal to achieve the optimal control and optimal state, based on the discretized dual problem, the primal problem is analyzed, in which a new discretized scheme for the $L^1$ norm is presented.
Compared the new discretized scheme $\|u_h\|_{\widetilde{L}_h^1(\Omega_h)}$ for $L^1$ norm with the nodal quadrature formula $\|u_h\|_{L_h^1(\Omega_h)}$, it is obvious that utilizing $\|u_h\|_{\widetilde{L}_h^1(\Omega_h)}$ to approximate $\|u_h\|_{L^1(\Omega_h)}$ is better than using $\|u_h\|_{L_h^1(\Omega_h)}$ in term of the order of approximation. Finally, finite element error estimates results for the primal problem with the new discretized scheme $\|u_h\|_{\widetilde{L}_h^1(\Omega_h)}$ are provided, which confirm that this approximation scheme will not change the order of error estimates.\\

\noindent {\bf Acknowledgments.} The authors would like to thank Prof. Defeng Sun at The Hong Kong Polytechnic University and Prof. Kim-Chuan Toh at National University of Singapore for their valuable suggestions that led to improvement in this paper. The research of the first author was supported by National Natural Science Foundation of China under Grant No. 91230103, 11571061, 11401075.

\end{document}